\documentclass[notitlepage]{amsart}
\usepackage{amsmath}%
\usepackage{amsfonts}%
\usepackage{amssymb}%
\usepackage{amsthm}
\usepackage{stmaryrd}
\usepackage[pdftex]{graphicx,color}
\usepackage{hyperref}
\usepackage{cite}
\usepackage[all]{xy}
 \newtheorem{theorem}{Theorem}
 \newtheorem{lemma}{Lemma}
 \newtheorem{proposition}{Proposition}
 \theoremstyle{definition}
 \newtheorem{definition}{Definition}
 \theoremstyle{remark}
 \newtheorem{example}{Example}
 \newtheorem{remark}{Remark}

\newcommand{\lb}{\llbracket}
\newcommand{\rb}{\rrbracket}
\newcommand{\EE}{\mathbb{E}}
\newcommand{\A}{\mathbb{A}}
\newcommand{\TT}{\mathbb{T}}

\newcommand{\id}{\mathsf{id}}
\newcommand{\Lied}{\mathcal{L}}
\newcommand{\dd}{\mathfrak{d}}

\newcommand{\f}{\mathfrak{f}}
\newcommand{\g}{\mathfrak{g}}
\newcommand{\h}{\mathfrak{h}}
\newcommand{\mf}{\mathfrak}
\newcommand{\mc}{\mathcal}
\newcommand{\ad}{\mathbf{ad}}

\newcommand{\on}{\operatorname}
\newcommand{\DIF}{\mathcal{D}}

\title[Quasi-Hamiltonian Groupoids]{Quasi-Hamiltonian Groupoids and Multiplicative Manin Pairs}

\author{David Li-Bland}
\author{Pavol \v{S}evera}

\address{Department of Mathematics, University of Toronto, 40 St George Street, Toronto, Ontario
M4S2E4, Canada}
\email{dbland@math.toronto.edu}
\address{Department of Mathematics, Universit\'{e} de Gen\`{e}ve, Geneva, Switzerland, on leave
from Dept. of Theoretical Physics, FMFI UK, Bratislava, Slovakia}
\email{pavol.severa@gmail.comphs}

\begin{document}
\begin{abstract}
We reformulate notions from the theory of quasi-Poisson $\g$-manifolds in terms of graded Poisson geometry and graded Poisson-Lie groups and prove that quasi-Poisson $\g$-manifolds integrate to quasi-Hamiltonian $\g$-groupoids.
We then interpret this result within the theory of  Dirac morphisms and multiplicative Manin pairs, to connect our work with more traditional approaches, and also to put it into a wider context suggesting possible generalizations.
\end{abstract}
\maketitle
\tableofcontents
\section*{Introduction}Let $G$ denote a Lie group whose Lie algebra $\g$ is equipped with an invariant inner product. Quasi-Hamiltonian $G$-manifolds were introduced in \cite{Alekseev97}, where they were shown to be equivalent to the theory of infinite dimensional Hamiltonian loop group spaces. In particular, they were used to simplify the study of the symplectic structure on the moduli space of flat connections by using finite dimensional techniques. In \cite{Alekseev99}, the more general quasi-Poisson $\g$-manifolds were introduced. These notions were further generalized and studied in subsequent papers \cite{Alekseev00,Bursztyn03-1,Bursztyn03,Xu03,Ponte05,QuasiPoissonAs,Bursztyn07-1,PureSpinorsOnL,Ponte08,Bursztyn08} (to cite a few).

The main objective of this paper is to prove that the following facts hold for an arbitrary quasi-Poisson $\g$-manifold, $M$:
\begin{enumerate}
\renewcommand{\labelenumi}{qP-\arabic{enumi}}
\item $T^*M$ inherits a Lie algebroid structure.\label{qP1}
\item If the Lie algebroid $T^*M$ integrates to a Lie groupoid $\Gamma\rightrightarrows M$, then $\Gamma$ inherits a quasi-Hamiltonian $\g$-structure. Moreover, the source map $s:\Gamma\to M$ is a quasi-Poisson morphism, while the target map $t:\Gamma\to M$ is anti-quasi-Poisson.\label{qP3}
\end{enumerate}
 Besides completing the theory of quasi-Poisson manifolds, our result can provide a new angle towards integrating certain Poisson structures. In addition to this, it has applications towards the integration of certain Courant algebroids. We plan to explore these consequences in forthcoming papers.

The methods we use to prove the results (qP-\ref{qP1} and \ref{qP3}) are of independent interest. We reformulate notions from the theory of quasi-Poisson $\g$-manifolds in terms of graded Poisson geometry and graded Poisson-Lie groups. Then we prove the results (qP-\ref{qP1} and \ref{qP3}) using well known theorems established for Poisson manifolds and Poisson Lie groups.  In particular, we prove (qP-\ref{qP3}) by interpreting structures in terms of Lie algebroid/groupoid morphisms, as was done in \cite{Mackenzie-Xu94,(MultDiracSt,Bursztyn09,Ponte05}. As a result, we are able to avoid infinite dimensional path spaces. We hope our approach will provide the reader with a fresh and insightful perspective on the theory of quasi-Poisson $\g$-manifolds.

The methods we use to derive the results (qP~\ref{qP1} and \ref{qP3}) are an application of the more general theory of MP-groupoids. We use the latter half of our paper to describe the theory of MP-groupoids. MP-groupoids are a reinterpretation of multiplicative Manin pairs \cite{(MultDiracSt} in terms of graded Poisson geometry \cite{NonComDiffForm,Bursztyn08}. In this portion of the paper, we also describe the link between our approach to quasi-Poisson $\g$-manifolds and the standard approach in terms of Dirac morphisms and multiplicative Manin pairs \cite{Bursztyn03,Bursztyn03-1,PureSpinorsOnL,Bursztyn08,(MultDiracSt}.

\subsection*{Overview}
Our paper is organized as follows. In \S~\ref{sect1}, we briefly summarize some background material and introduce our main result, the integration of quasi-Poisson $\g$-manifolds. In \S~\ref{part1}, we provide a new perspective on the theory of quasi-Poisson $\g$-manifolds using the theory graded Poisson geometry and graded Poisson-Lie groups. We then use this new perspective to prove the results from \S~\ref{sect1}. Next, in \S~\ref{examples1}, we provide some detailed examples of the integration of quasi-Poisson $\g$-manifolds.

The remainder of the paper is spent relating our approach to the theory of Manin pairs. In \S~\ref{sect:CAandMp} we review the definitions of Courant algebroids, Manin pairs, their morphisms and the category of multiplicative Manin pairs. Following this, in \S~\ref{sect:MpandMPmflds} we recall the relationship between the categories of Manin pairs and graded Poisson manifolds. Using this relationship, we introduce the infinitesimal notion corresponding to a multiplicative Manin pair. We apply these concepts in \S~\ref{sect:MpandqPstruct} to relate the content of \S~\ref{part1} to the theory of Manin pairs. 

\subsection*{Acknowledgements}
We would like to thank Eckhard Meinrenken for his advice and suggestions. D.L.-B. was supported by an NSERC CGS-D Grant, and thanks the Universit\'{e} de Gen\`{e}ve, and Anton Alekseev in particular, for their hospitality during his visit. P.~\v S. was supported by the Swiss National Science Foundation (grant 200020-120042 and 200020-126817). Finally, we would like to thank the referees for their helpful comments.

\section{Background and statement of results}\label{sect1}
In this section we want to recall the theory of quasi-Poisson $\g$-manifolds. To provide some intuition and motivate the definitions, we will develop both the theory of quasi-Poisson $\g$-manifolds and the theory of Poisson manifolds in parallel. We will describe the Poisson case in a series of remarks.

In \S~\ref{part1} we will use a graded version of the theory of Poisson manifolds and Lie bialgebra actions to prove some results about quasi-Poisson $\g$-manifolds. In particular, for \S~\ref{part1} we will require some understanding of quasi-triangular Lie bialgebras. We intend to summarize this background material in this section.

\subsection{Quasi-Poisson $\g$-manifolds}

Let $\g$ be a Lie algebra with a chosen ad-invariant element $s\in S^2\g$. Define $s^\sharp:\g^*\to \g$ by $\beta(s^\sharp(\alpha))=s(\alpha,\beta)$ for $\alpha,\beta\in\g^*$. Let $\phi\in\wedge^3\g$ be given by 
$$\phi(\alpha,\beta,\gamma)=\frac{1}{2}\alpha([s^\sharp \beta, s^\sharp \gamma])\quad (\alpha,\beta,\gamma\in\g^*).$$

\begin{definition}\cite{Alekseev99, Alekseev00} A quasi-Poisson $\g$-manifold is triple $(M,\rho,\pi)$, where
\begin{itemize}
\item $M$ is a $\g$-manifold,
\item $\rho:M\times\g\to TM$ is the anchor map for the action Lie algebroid, and
\item  $\pi\in \mathbf{\Gamma}(\wedge^2 TM)^\g$ is a $\g$-invariant bivector field
\end{itemize} satisfying
  \begin{equation}\label{qPoissonCond}[\pi,\pi]=\rho(\phi),\quad\text{ and }\quad [\pi,\rho(\xi)]=0,\end{equation} for any $\xi\in\g$.
Here we view $\phi$ and $\xi$ as constant sections of $M\times\wedge\g$ and extend $\rho$ to a morphism $\rho:M\times\wedge\g\to \wedge TM$. The bracket in \eqref{qPoissonCond} is the Schouten bracket on multivector fields.


Let $(M_i,\rho_i,\pi_i)$ be two quasi-Poisson $\g$-manifolds (for $i=1,2$). A map $f:M_1\to M_2$ is a quasi-Poisson morphism if $\pi_2=f_*\pi_1$ and $\rho_2=f_*\circ\rho_1$.
\end{definition}

\begin{example}\label{coisotropicqpm}Suppose that $\rho$ is an action of $\g$ on $M$. For $x\in M$ let $\g_x=\ker(\rho_x)$ be the stabilizer at the point $x$, and $\g_x^\perp\subset\g^*$ be its annihilator. Then $\rho(\phi)=0$ if $s^\sharp(\g_x^\perp)\subset \g_x$, in which case we say the stabilizers are coisotropic (with respect to $s$). In this case, the triple $(M,\rho,0)$ is a quasi-Poisson $\g$-manifold.  
\end{example}

\begin{remark}[Poisson Parallel]\label{rem:ManTrip}We now recall the related notion in Poisson geometry.

Let $\h$ be a Lie bialgebra with cobracket $\delta:\h\to\wedge^2\h$. A Poisson $\h$-manifold is a triple $(N,\rho,\pi)$, where
\begin{itemize}
\item $N$ is a $\h$-manifold,
\item $\rho:M\times\h\to TN$ is the anchor map for the action Lie algebroid, and
\item $\pi\in\mathbf{\Gamma}(\wedge^2TN)$ is a bivector field,
\end{itemize} satisfying  \begin{equation}\label{PoissonCond}[\pi,\pi]=0,\quad\text{ and }\quad[\pi,\rho(\xi)]=\rho(\delta(\xi))\end{equation}
 for every $\xi\in \h$, which we view as a constant section of $M\times\h$. 

If $H$ denotes the Poisson-Lie group corresponding to the Lie bialgebra $\h$, and if $\rho$ can be integrated to an action of $H$ on $N$ (or we work with local actions), the action $H\times N\to N$ becomes a Poisson map \cite{thesis-3,lu90}.



When $\h={0}$, then we may refer to a Poisson $\h$-manifold $(N,\rho=0,\pi)$ as a \emph{Poisson manifold} $(N,\pi)$.
\end{remark}


\begin{remark}\label{qaction_vs_bialg}
Recall that the Lie bracket on $\g$ extends to a Gerstenhaber bracket on $\bigwedge\g$. Suppose there is an element $u\in\bigwedge^2 \g$ such that $[u,u]=-\phi$.
Then $\g$ becomes a Lie bialgebra with cobracket $\delta:\xi\to[u,\xi]$. Lie bialgebras of this type are called {\em quasi-triangular} \cite{KS92,KS93} and the combination $r=s+u\in \g\otimes\g$ is called a {\em classical $r$-matrix}.
 
In this case, quasi-Poisson $\g$-manifolds and Poisson $\g$-manifolds are equivalent via a ``twist'' by $u$, as shown in \cite{Alekseev99}. To recall the details, if $(M,\rho,\pi)$ is a quasi-Poisson $\g$-manifold then $\pi'=\pi+\rho(u)$ satisfies $[\pi',\pi']=0$, i.e.\ $\pi'$ is a Poisson structure. Moreover the action of $\g$ via $\rho$ becomes an action of the Lie bialgebra $\g$ on the Poisson manifold $(M,\pi')$.
  
Note that in \S~\ref{part1} we will use a graded version of this equivalence.
\end{remark}

For simplicity, we shall {\em restrict from now on} to the case where {\em $s$ is non-degenerate}. We shall identify $\g$ with $\g^*$ via $s$, and let $\langle\cdot,\cdot\rangle$ denote the corresponding $\ad$-invariant inner-product on $\g$. Such a Lie algebra is called \emph{quadratic}. Let $e_i$ and $e^i$ denote two bases of $\g$ dual with respect to $\langle\cdot,\cdot\rangle$.

\begin{remark}
Suppose $(M,\rho,0)$ is as in Example~\ref{coisotropicqpm}, $s$ is non-degenerate, and $\g=\f\oplus\h$ where $\f$ and $\h$ are Lagrangian subalgebras. We can choose bases $\{f_i\}$ of $\f$ and $\{h_i\}$ of $\h$ such that $\langle f_i,h_i\rangle=1$ and $\langle f_i,h_j\rangle=0$ for $i\neq j$. Then with $u:=\sum_i f_i\wedge h_i$, $r=s+u\in \g\otimes\g$ is a classical $r$-matrix. As in Remark~\ref{qaction_vs_bialg}, $\pi':=\pi+\rho(u)$ defines a Poisson structure on $M$.

The Poisson structure $\pi'$ was constructed in \cite{Li-Bland08} via a Courant algebroid structure on $M\times\g$. The approach via quasi-Poisson $\g$-manifolds appears to be more direct.
\end{remark}

Recall that Lie algebroid structure on a vector bundle $A\to M$ is equivalent to a differential on the algebra $\mathbf{\Gamma}(\wedge^* A^*)$ \cite{LieAlgebroidsH,Courant90}.

\begin{theorem}\label{liealgstruct}
If $(M,\rho,\pi)$ is a quasi-Poisson $\g$-manifold, then $T^*M$ becomes a Lie algebroid, where the Lie algebroid differential $d_{T^*M}$ on $\mathbf{\Gamma}(\wedge^* TM)$ is
\begin{equation}\label{liealgdiff}d_{T^*M}=[\pi,\cdot]+\frac{1}{2}\sum_i\rho(e^i)\wedge[\rho(e_i),\cdot]\end{equation}
Furthermore, the induced action of $\g$ on $T^*M$ preserves the Lie algebroid structure.

Let $\on{pr}_\g:M\times\g\to\g$ be the projection to the second factor, and $\rho^*:T^*M\to M\times\g^*\cong M\times \g$ the transpose of $\rho$, then $\mu_\rho:=\on{pr}_\g\circ\rho^*:T^*M\to \g$ is a Lie algebroid morphism.
\end{theorem}
 This can be proven by a direct calculation: the two parts of $d_{T^*M}$ commute with each other and their squares cancel each other. We give a conceptual proof of the theorem in \S~\ref{mathbbR2timesg1}.

The corresponding Lie bracket on 1-forms $\alpha,\beta\in\Omega^1(M)$ is 
$$[\alpha,\beta]=[\alpha,\beta]_\pi+\frac{1}{2}\sum_i\bigg(\alpha\big(\rho(e^i)\big)\Lied_{\rho(e_i)}\beta-\beta\big(\rho(e^i)\big)\Lied_{\rho(e_i)}\alpha\bigg),$$
where \begin{equation}\label{brackpi}[\alpha,\beta]_\pi=d\pi(\alpha,\beta)+\iota_{\pi^\sharp(\alpha)}d\beta-\iota_{\pi^\sharp(\beta)}d\alpha\end{equation} is the Koszul bracket (here $\beta(\pi^\sharp(\alpha))=\pi(\alpha,\beta)$).
The anchor map, $\mathbf{a}:T^*M\to TM$, is 
\begin{equation}\label{anchor}\mathbf{a}=\pi^\sharp+\frac{1}{2}\rho\circ\rho^*.\end{equation}


We call a quasi-Poisson $\g$-manifold \emph{integrable} if the Lie algebroid structure on the cotangent bundle is integrable to a Lie groupoid.

\begin{remark}[Poisson Parallel]

If $(N,\pi)$ is a Poisson manifold, there is a Lie algebroid structure on $T^*N$ whose corresponding Lie algebroid differential $d_{T^*N}:\mathbf{\Gamma}(\wedge^n TN)\to\mathbf{\Gamma}(\wedge^{n+1}TN)$ is $$d_{T^*N}=[\pi,\cdot].$$
In this case, the anchor map $\mathbf{a}:T^*M\to TM$, is given by \begin{equation}\label{panchor}\mathbf{a}=\pi^\sharp,\end{equation} and the bracket is given by \eqref{brackpi}.

If a Lie bialgebra $\h$ acts on $N$, then $\mu_\rho:=\on{pr}_{\h^*}\circ\rho^*:T^*N\to \h^*$ is a Lie algebroid morphism, where $\on{pr}_{\h^*}:M\times \h^*\to \h^*$ is the projection to the second factor and $\rho^*:T^*M\to M\times\h^*$ is the transpose of $\rho$. 
\end{remark}

\begin{remark}\label{LieAlgL}
It was shown in \cite{Bursztyn03} that whenever $(M,\rho,\pi)$ is a quasi-Poisson $\g$-manifold, there is a Lie algebroid structure on $A=T^*M\oplus\g$. The Lie algebroid $T^*M$ described in Theorem~\ref{liealgstruct} is embedded as a subalgebroid of $T^*M\oplus\g$ by the map $\alpha\to\alpha+\rho^*(\alpha)$.

The following was pointed out by a referee: Suppose that $(M,\rho,\pi)$ also possesses a moment map $\Phi:M\to G$ (see Definition~\ref{MomentMapDefn}). Let $\theta$ denote the (left) Maurer-Cartan form on $G$, and let $\eta=\frac{1}{2}\langle [\theta,\theta],\theta\rangle$. In \cite{Bursztyn03}, H. Bursztyn and M. Crainic describe a $\Phi^*\eta$-twisted Dirac structure $L\subset TM\oplus T^*M$ associated to the quasi-Poisson $\g$-structure on $M$. In \cite[Proposition 3.19]{Bursztyn03} they describe a map $$A\to L\subset TM\oplus T^*M,$$ which restricts to $T^*M\cong\on{Gr}_{\rho^*}\subset A$ to define an isomorphism of Lie algebroids $T^*M\cong L.$

In particular, \cite{Bursztyn03} shows that in the presence of a moment map the Lie algebroid described in Theorem~\ref{liealgstruct} can be viewed as a $\Phi^*\eta$-twisted Dirac structure.
\end{remark}

\begin{remark}
The foliation of the quasi-Poisson $\g$-manifold $(M,\rho,\pi)$ given by the Lie algebroid $T^*M$ is \emph{different} from the foliation given in \cite{Alekseev99,Alekseev00,Bursztyn03}. The latter foliation is tangent to $\rho_x(\g)+\pi_x^\sharp T^*_xM$ at any point $x\in M$; in particular, the leaves contain the $\g$-orbits. This is not the case for the foliation given by the Lie algebroid $T^*M$; for instance, if as in Example~\ref{coisotropicqpm}, $(M,\rho,0)$ is a quasi-Poisson $\g$-manifold, the anchor map $\mathbf{a}:T^*M\to TM$ is trivial ($\rho_x\circ\rho^*_x=0$ for any point $x\in M$, since the stabilizers of $\rho$ are coisotropic), while the $\g$-orbits may not be. 

On the other hand, as we shall see below, for a Hamiltonian quasi-Poisson $\g$-manifold these two foliations coincide.
\end{remark}

\begin{definition}\label{qsympl}
A \emph{quasi-symplectic} $\g$-manifold is a quasi-Poisson $\g$-manifold $(M,\pi,\rho)$ such that the anchor map \eqref{anchor} is bijective.
\end{definition}

\begin{remark}[Poisson Parallel]
A Poisson manifold $(N,\pi)$ is called \emph{symplectic} if the corresponding anchor map \eqref{panchor} is bijective.
\end{remark}

\subsection{Hamiltonian quasi-Poisson $\g$-manifolds}
There is a concept of a group valued moment map for quasi-Poisson $\g$-manifolds \cite{Alekseev99,Alekseev00}. Let $G$ be a Lie group with Lie algebra $\g$. For any $\xi\in\g$ let $\xi^L,\xi^R\in\mathbf{\Gamma}(TG)$ denote the corresponding left and right invariant vector fields. Let $\theta^L,\theta^R\in\Omega^1(G,\g)$ denote the left and right Maurer Cartan forms on $G$ defined by
$$\theta^L(\xi^L)=\xi,\quad \theta^R(\xi^R)=\xi.$$


\begin{definition}\label{MomentMapDefn}\cite{Alekseev99,Alekseev00} A map $\Phi:M\to G$ is called a \emph{moment map} for the quasi-Poisson $\g$-manifold $(M,\rho,\pi)$ if
\begin{itemize}
\item $\Phi$ is $\g$-equivariant, and 
\item $\pi^\sharp(\Phi^*(\alpha))=\rho(\Phi^*(b^*\alpha))\text{ for any }\alpha\in\Omega^1(G)$,
\end{itemize}
 where the vector bundle map $b:G\times\g\to TG$ is given by \begin{equation}\label{beq}b:(g,\xi)\to \frac{1}{2}\big(\xi^L(g)+\xi^R(g)\big).\end{equation}

Under these conditions, we call the quadruple $(M,\rho,\pi,\Phi)$ a \emph{Hamiltonian quasi-Poisson $\g$-manifold}, or a Hamiltonian quasi-Poisson $\g$-structure on $M$.
\end{definition}


\begin{theorem}\label{g-inclusion}
If $(M,\rho,\pi,\Phi)$ is a Hamiltonian quasi-Poisson $\g$-manifold, then the map
\begin{equation}\label{i-map-def}i:\g\to \Omega^1(M),\quad i(\xi)=\Phi^*\langle\xi,\theta^L\rangle\end{equation}
 is a morphism of Lie algebras such that $\mathbf{a}\circ i=\rho$.
\end{theorem}
Again this can be proved by a direct calculation, and we give a conceptual proof in \S~\ref{qhqpgmr}.

\begin{remark}[Poisson Parallel]
In \S~\ref{part1}, we will need the corresponding notion of group valued moment maps for Poisson geometry \cite{thesis-3}, summarized as follows.

Let again $\h$ be a Lie bialgebra, and let $H^*$ denote the 1-connected Poisson Lie group integrating $\h^*$. Recall \cite{thesis-3} that a Poisson map
$$\Phi:N\to H^*$$
gives rise to a morphism of Lie algebras
$$i:\h\to\Omega^1(N),\quad i(\xi)=\Phi^*(\xi_L)$$
where $\xi_L$ is the left-invariant 1-form on $H^*$ equal to $\xi$ at the group unit. The morphism $i$ then in turn produces an action $\rho$ of $\h$ on $N$ via the anchor map on $T^*N$, i.e.
$$\rho(\xi)=\pi^\sharp(i(\xi)).$$
The map $\Phi$ is automatically $\h$-equivariant, where the so-called dressing action of $\h$ on $H^*$ comes from the identity moment map $H^*\to H^*$. $\Phi$ is called a moment map for the action $\rho$.
\end{remark}

\begin{remark}\label{rem:mmdeg}
When $s\in (S^2\g)^\g$ is not assumed to be non-degenerate, one can proceed as follows. The element $s$ is equivalent to a triple $(\mf{d},\mf{g},\mf{g}')$, where \begin{itemize}

\item $\mf{d}$ is a quadratic Lie algebra,
\item $\g\subset\dd$ is a Lagrangian subalgebra,
\item $\g'\subset\dd$ an ideal such that $\dd=\g\oplus\g'$ as vector spaces, and
\end{itemize}
the restriction to $\g'$ of the inner product on $\mf{d}$ is $s$ (where we identify $\g'$ with $\g^*$ via the inner product in $\dd$).
          
Moment maps then have value in a group $G'$ integrating $\g'$.
\end{remark}




\begin{theorem}\label{thm:crit_for_qH}
 If $(M,\rho,\pi,\Phi)$ is a Hamiltonian quasi-Poisson $\g$-manifold then
 $$\rho_x(\g)+ \pi_x^\sharp(T^*_xM)=\mathbf{a}(T^*_xM).$$ 
\end{theorem}

\begin{proof}
Recall that $\mathbf{a}=\pi^\sharp+\frac{1}{2}\rho\circ\rho^*$, hence $\rho_x(\g)+ \pi_x^\sharp(T^*_xM)\supseteq \mathbf{a}(T^*_xM)$. On the other hand, by Theorem~\ref{g-inclusion} we have $\mathbf{a}\circ i=\rho$. Hence, $$\pi^\sharp=\mathbf{a}-\frac{1}{2}\rho\circ\rho^*=\mathbf{a}\circ(\id-\frac{1}{2} i\circ\rho^*),$$
and $\pi_x^\sharp(T^*_xM)\subseteq\mathbf{a}_x(T^*_xM)$. 

\end{proof}

\begin{remark}
 A Hamiltonian quasi-Poisson $\g$-manifold $(M,\rho,\pi,\Phi)$  is a \emph{quasi-Hamiltonian $\g$-manifold} \cite{Alekseev97,Alekseev99,Alekseev00,PureSpinorsOnL,Bursztyn03}, if for every point $x\in M$,
\begin{equation}\label{nondegenqh}\rho_x(\g)+ \pi_x^\sharp(T^*_xM)=T_xM.\end{equation}



It follows from Theorem~\ref{thm:crit_for_qH} that a quasi-Hamiltonian $\g$-manifold is equivalent to a Hamiltonian quasi-symplectic $\g$-manifold. We will use the latter term in this paper.

\end{remark}

\subsection{Fusion}\label{FusionSect}

  The category of quasi-Poisson $\g$-manifolds, has a braided monoidal structure given by fusion \cite{Alekseev00}. Let $(M,\rho,\pi)$ be a quasi-Poisson $\g\oplus\g$-manifold, \begin{equation}\label{psi}\psi=\frac{1}{2}\sum_i(e^i,0)\wedge(0,e_i)\in \wedge^2(\g\oplus\g),\end{equation} and let $\on{diag}(\g)\cong \g$ denote the diagonal subalgebra of $\g\oplus\g$. The quasi-Poisson $\g$-manifold
 \begin{equation}\label{fusion}(M,\rho\rvert_{\on{diag}(\g)},\pi_{\on{fusion}}),\quad \text{ with }\pi_{\on{fusion}}=\pi+\rho(\psi),\end{equation} is called the fusion of $M$ \cite{Alekseev00}. Moreover, if $(\Phi_1,\Phi_2):M\to G\times G$ is a moment map, then the pointwise product $\Phi_1\Phi_2:M\to G$ is a moment map for the fusion \cite{Alekseev00}.

If $(M_i,\rho_i,\pi_i,\Phi_i)$ (for $i=1,2$) are two Hamiltonian quasi-Poisson $\g$-manifolds, then $(M_1\times M_2,\rho_1\times\rho_2,\pi_1+\pi_2,\Phi_1\times\Phi_2)$ is a quasi-Poisson $\g\oplus\g$-manifold. The fusion $$(M_1\times M_2,(\rho_1\times\rho_2)\rvert_{\on{diag}(\g)},(\pi_1+\pi_2)_{\on{fusion}},\Phi_1\Phi_2)$$ is called the fusion product, and denoted $$(M_1,\rho_1,\pi_1)\circledast(M_2,\rho_2,\pi_2),$$ or just $M_1\circledast M_2$.

The fusion product of two quasi-Poisson $\g$-manifolds is defined similarly, one just ignores the moment maps.

The monoidal category of Hamiltonian quasi-Poisson $G$-manifolds is braided \cite{Alekseev00}: if $(M_i,\rho_i,\pi_i,\Phi_i)$ are two Hamiltonian quasi-Poisson $G$-manifolds, the corresponding isomorphism between fusion products
$$M_1\circledast M_2 \to M_2 \circledast M_1$$
is given by $(x_1,x_2)\mapsto (\Phi_1(x_1)\cdot x_2, x_1)$.

To provide an alternate explanation for this monoidal structure in \S~\ref{part1}, we will need to understand the story for Poisson manifolds.

\begin{remark}[Poisson Parallel]
If $(M,\pi_M)$ and $(N,\pi_N)$ are two Poisson manifolds, then $(M\times N,\pi_M+\pi_N)$ is also a Poisson manifold. If $M\to H^*$, $N\to H^*$ are Poisson (moment) maps, we can compose
$$M\times N\to H^*\times H^*\to H^*$$
(the latter arrow is the product in $H^*$)
 to get a Poisson map $M\times N\to H^*$. The category of those Poisson manifolds with a Poisson map to $H^*$ is thus monoidal (but not necessarily braided). Notice, that unlike the case of quasi-Poisson manifolds, the resulting action of $\h$ on $M\times N$ is not just the diagonal action -- it is twisted by the moment map on $M$. On the other hand, the Poisson bivector is simply the sum $\pi_{M_1}+\pi_{M_2}$.

If $(M_1,\pi_1)$ and $(M_2,\pi_2)$ are Poisson manifolds, then the Lie algebroid structure on $T^*(M_1\times M_2)$ is the direct sum of the Lie algebroids $T^*M_1$ and $T^*M_2$. 
\end{remark}

For quasi-Poisson $\g$-manifolds, one may ask how the Lie algebroids $T^*(M_1\circledast M_2)$ and $T^*M_1\oplus T^*M_2$ are related. A direct computation shows
\[d_{T^*(M_1\circledast M_2)}=d_{T^*M_1}+d_{T^*M_2}+\sum_i\rho_1(e^i)\wedge[\rho_2(e_i),\cdot]\ ,\]
so that the obvious isomorphism of vector spaces $T^*(M_1\circledast M_2)\cong T^*M_1\oplus T^*M_2$ is not an isomorphism of Lie algebroids. However, for \emph{Hamiltonian} quasi-Poisson manifolds, there is a non-standard isomorphism $T^*(M_1\times M_2)\cong T^*M_1\oplus T^*M_2$ which is an isomorphism of Lie algebroids.

By a \emph{comorphism} \cite{Mackenzie05} from a Lie algebroid $A\to M$ to a Lie algebroid $A'\to M'$ we mean a morphism of Gerstenhaber algebras
$\mathbf{\Gamma}(\bigwedge A')\to\mathbf{\Gamma}(\bigwedge A)$. The assignment $M\mapsto T^*M$ is a functor from the category of quasi-Poisson $\g$-manifolds and quasi-Poisson morphisms to the category of Lie algebroids and comorphisms.

\begin{proposition}\label{prop:lalg-mon}
Let $(M_i,\rho_i,\pi_i,\Phi_i)$ $(i=1,2)$ be two Hamiltonian quasi-Poisson $\g$-manifolds and let the isomorphism $J:T^*M_1\oplus T^*M_2\to T^*(M_1\circledast M_2)$ be given by $(\alpha,\beta)\mapsto(\alpha,\beta-i_2(\rho_1^*(\alpha)))$. Then $J$ is a isomorphism of Lie algebroids.
Moreover, $J$ is a natural transformation which makes the functor $M\mapsto T^*M$ (from the category of Hamiltonian quasi-Poisson $\g$-manifolds to the category of Lie algebroids and comorphisms) strongly monoidal.
\end{proposition}
The proof is in \S~\ref{sect:ghat-qtr}, but it requires a deeper understanding of the relationship between the monoidal structures for the categories of quasi-triangular Poisson manifolds and quasi-Poisson $\g$-manifolds, which we shall now recall.

Let $\h$ be a quasi-triangular Lie-bialgebra. As described in Remark~\ref{rem:ManTrip}, $\h$ corresponds to a Manin triple $(\dd,\h,\h^*)$. Let $\h'\subset\dd=\h\oplus\h^*$ be the graph of the $r$-matrix. Then $\h'$ is an ideal, so that $(\dd,\h,\h')$ is as in Remark \ref{rem:mmdeg}. Let $H,H^*,H'\subset D$ be groups with Lie algebras $\h,\h^*,\h'\subset\dd$. Suppose further that the maps $H^*\to D/H$ and $H'\to D/H$ are bijections. 

Suppose $(M,\rho,\pi,\Phi)$ is a quasi-Poisson $H$-manifold with $H$-action $\rho:H\times M\to M$, bivector $\pi$, and with moment map $\Phi:M\to D/G\cong H'$. Define $F_\h(M,\rho,\pi,\Phi):=(M,\rho,\pi',\Phi)$ to be the Poisson $H$-manifold with bivector $\pi'$ given by twisting $\pi$ as in in Remark~\ref{qaction_vs_bialg}, the same $H$-action, and the same moment map $\Phi:M\to D/G\cong H^*$.  Then as shown in \cite{Alekseev99,Bursztyn07-1,Ponte08,Bursztyn08}, the functor $$F_\h:\on{Ham-qPois}_\h\to\on{Ham-Pois}_\h$$ describes an equivalence between the category $\on{Ham-qPois}_\h$ of quasi-Poisson $H$-manifolds with $H'$-valued moment maps and the category $\on{Ham-Pois}_\h$ of Poisson $H$-manifolds with $H^*$-valued moment maps. 

There is a canonical choice for the inverse functor, $F^{-1}_\h(M,\rho,\pi',\Phi):=(M,\rho,\pi,\Phi)$, where $\pi$ is constructed from $\pi'$ by reversing the procedure in Remark~\ref{qaction_vs_bialg}.

\begin{remark}
This equivalence can be understood more intrinsically. What is significant is that there is a natural Dirac structure living over $D/G$,   $$(\dd\times D/G,\g\times D/G)$$ (see Example~\ref{ex:s-val} for more details). It was shown in \cite{Bursztyn07-1,Ponte08,Bursztyn08} that a morphism of Manin Pairs $$(\Phi,K):(\TT M,TM)\dasharrow (\dd\times D/G,\g\times D/G)$$ is the intrinsic data underlying both $(M,\rho,\pi,\Phi)$ and $(M,\rho,\pi',\Phi)$. Indeed, the action $\rho$ is specified by the morphism of Manin pairs, and the bivectors $\pi$ and $\pi'$ arise from choosing two different Lagrangian complements to $\g$ in $\dd$ \cite{Alekseev99,Ponte08}.
\end{remark}




\begin{remark}\label{rem:mmequiv}
As shown by A.~Weinstein and P.~Xu \cite{weinxu92}, when the Lie-bialgebra $\h$ is quasi-triangular, the category of Poisson $H$-manifolds with $H^*$-valued moment maps is braided monoidal. In fact, $F_\h$ is a strong monoidal functor with the natural transformation given by
\[\mathcal{J}:F_\h(M_1\circledast M_2)\to F_\h(M_1)\times F_\h(M_2),\quad (x_1,x_2)\mapsto(x_1,\,j(\Phi_1(x_1))\cdot x_2)\ ,\]
where $\Phi_1:M_1\to D/G\cong H^*$ is the moment map for $M_1$ and $j:H^*\to H$ is the map specified by the condition 
\[g\, j(g)\in H'\text{ for every }g\in H^*.\]
(To define $j:H^*\to H$, we used the fact that the maps $H^*\to D/H$ and $H'\to D/H$ are bijections.)


We shall need a graded version of this fact in the proof of Proposition~\ref{prop:lalg-mon} in \S~\ref{sect:ghat-qtr}.


\end{remark}

\subsection{Hamiltonian quasi-Poisson $\g$-groupoids}\label{qhqpg1}

Let $\Gamma\rightrightarrows M$ be a group\-oid, and let $\on{Gr}_{mult_\Gamma}=\{(g,h,g\cdot h)\}\subset \Gamma\times\Gamma\times\Gamma$ denote the graph of the multiplication map. A bivector field  $\pi_\Gamma\in\mathbf{\Gamma}(\wedge^2  T\Gamma)$ is said to be \emph{multiplicative} \cite{Ponte05} if $\on{Gr}_{mult_\Gamma}$ is a coisotropic submanifold of $(\Gamma,\pi_\Gamma)\times(\Gamma,\pi_\Gamma)\times(\Gamma,-\pi_\Gamma)$.




\begin{definition}\label{qhamg}
Suppose that $\Gamma\rightrightarrows M$ is a groupoid, and $(\Gamma,\rho,\pi_\Gamma,\Phi)$ is a Hamiltonian quasi-Poisson $\g$-manifold. It is called a \emph{Hamiltonian quasi-Poisson $\g$-groupoid} if
\begin{itemize}
\item $\Phi:\Gamma\to G$ is a morphism of groupoids,
\item $\g$ acts on $\Gamma$ by (infinitesimal) groupoid automorphisms, and
\item $\on{Gr}_{mult_\Gamma}$ is coisotropic with respect to the bivector field $$((\pi_\Gamma+\pi_\Gamma)_{\on{fusion}})_{1,2}-(\pi_\Gamma)_3,$$ where $((\pi_\Gamma+\pi_\Gamma)_{\on{fusion}})_{1,2}$ appears on the first two factors of $\Gamma\times\Gamma\times\Gamma$ and $(\pi_\Gamma)_3$ appears on the third.

\end{itemize}
We refer to the last condition as $\pi$ being \emph{fusion multiplicative}.

A \emph{Hamiltonian quasi-symplectic $\g$-groupoid} is a Hamiltonian quasi-Poisson $\g$-groupoid such that the anchor map \eqref{anchor} is bijective.

A Hamiltonian quasi-Poisson $\g$-groupoid is called \emph{source 1-connected} if $\Gamma$ is source 1-connected and $G$ is 1-connected.

\end{definition}


\begin{remark}[Poisson Parallel]
If $(N,\pi)$ is a Poisson manifold, and the Lie algebroid $T^*N$ integrates to a (possibly local) Lie groupoid $\Gamma\rightrightarrows N$, then \cite{weinstein87-1}
\begin{itemize}
\renewcommand{\labelenumi}{P-\arabic{enumi}}
\item there is a bivector field $\pi_{\Gamma}\in\mathbf{\Gamma}(\wedge^2T\Gamma)$ such that $(\Gamma,\pi_{\Gamma})$ is a Poisson manifold,
\item $\pi_\Gamma$ is non-degenerate (so that $\Gamma$ is in fact symplectic).
\item $\pi_{\Gamma}$ is multiplicative, so that $(\Gamma,\pi_{\Gamma})$ is a Poisson groupoid \cite{weinstein87} (in fact, a symplectic groupoid).
\end{itemize}
Suppose, in addition, that a Lie bialgebra $\h$ acts on $N$. We can interpret the action as a Lie bialgebroid morphism $T^*N\to\h^*$, which then integrates to a Poisson groupoid morphism $\Gamma\to H^*$ (see \cite{Xu95}).
\end{remark}

\subsection{Main results}
We may now state the first of our main results

\begin{theorem}\label{mainthm1}
There is a one-to-one correspondence between source 1-connected Hamiltonian quasi-symplectic $\g$-groupoids $(\Gamma,\rho_\Gamma,\pi_\Gamma,\Phi)$ and integrable quasi-Poisson $\g$-manifolds $(M,\rho,\pi)$. Under this correspondence, the Lie algebroid of $\Gamma$ is $T^*M$ and $\Phi$ integrates the Lie algebroid morphism $\mu_\rho:T^*M\to\g$. Furthermore, the source map $s:\Gamma\to M$ is a quasi-Poisson morphism, while the target map $t:\Gamma\to M$ is anti-quasi-Poisson.
\end{theorem}

\begin{remark}
Theorem~\ref{mainthm1} was already established for the case where the quasi-Poisson $\g$-manifold $(M,\rho,\pi)$ possess a moment map $\Phi:M\to G$. This fact was pointed out to us by a referee, and we explain it in Remark~\ref{establishedResult} after recalling some background.
\end{remark}



Theorem~\ref{mainthm1} prompts one to ask what a general Hamiltonian quasi-Poisson $\g$-groupoid corresponds to infinitesimally. To answer this, we extend the notion of a quasi-Poisson $\g$-manifold, by replacing the tangent bundle with an arbitrary Lie algebroid.

\begin{definition}[quasi-Poisson $\g$-bialgebroid]\label{def:qpgb}
A quasi-Poisson $\g$-bialgebroid over a $\g$-manifold $M$ is a triple $(A,\rho,\DIF)$ consisting of 
\begin{itemize}
\item a Lie algebroid $A\to M$,
\item a Lie algebroid morphism $\rho:M\times\g\to A$, where $M\times\g$ is the action Lie algebroid, and
\item a degree $+1$ derivation $\DIF$ of the Gerstenhaber algebra $\mathbf{\Gamma}(\bigwedge A)$,
\end{itemize}
 such that 
\begin{itemize}
\item $\DIF\rho(\xi)=0$ for any constant section $\xi\mathbf{\Gamma}(M\times\g)$, and
\item $\DIF^2=\frac{1}{2}[\rho(\phi),\cdot]\,$, where we view $\phi$ as a constant section of $M\times\g$.
\end{itemize}

Let $\on{pr}_\g:M\times\g\to\g$ denote the projection to the second factor, and $\rho^*:A^*\to M\times\g^*\cong M\times\g$ be the transpose of $\rho$. Define $\mu_\rho:A^*\to\g$ by $\mu_\rho:=\on{pr}_\g\circ\rho^*$.
\end{definition}

\begin{remark}
 Quasi-Poisson $\g$-bialgebroids are examples of Lie quasi-bialgebroids, a notion introduced in \cite{Roytenberg99,Roytenberg02} by D. Roytenberg (see also \cite{QuasiTwistedAl}). A Lie quasi-bialgebroid is a triple $(A,\DIF,\chi)$, where $A$ is a Lie algebroid, $\DIF$ is a degree $+1$ derivation of the Gerstenhaber algebra $\mathbf{\Gamma}(\wedge A)$, and $\chi\in\mathbf{\Gamma}(\wedge^3 A)$. They must satisfy the equations $\DIF^2=\frac{1}{2}[\chi,\cdot]$ and $\DIF\chi=0$.

Therefore, if $(A,\rho,\DIF)$ is a quasi-Poisson $\g$-bialgebroid, then $(A,\DIF,\rho(\phi))$ is a Lie quasi-bialgebroid.
\end{remark}

\begin{example}[quasi-Poisson $\g$-manifolds]\label{QPMisQPBA}
Suppose that $(M,\rho,\pi)$ is a quasi-Poisson $\g$-manifold. Let $$\DIF_\pi=[\pi,\cdot]_{\on{Schouten}}$$ be the derivation of $\mathbf{\Gamma}(\wedge^* TM)$ given by the Schouten bracket. Then $(TM,\rho,\DIF_\pi)$ is a quasi-Poisson $\g$-bialgebroid.
\end{example}

\begin{proposition}\label{inducedqpstruct}
Let $A\to M$ be a Lie algebroid with anchor map $\mathbf{a}_A:A\to TM$. A compatible quasi-Poisson $\g$-bialgebroid structure $(A,\rho,\DIF)$ defines a canonical quasi-Poisson $\g$-structure $(M,\mathbf{a}_A\circ\rho,\pi_\DIF)$ on $M$ via
$$\pi_\DIF^\sharp(df)=\mathbf{a}_A\circ \DIF f,$$
where we view the function $f\in C^\infty(M)$ as an element of $\mathbf{\Gamma}(\wedge^0A)$.
\end{proposition}
The result is just a special case of \cite[Proposition 4.8]{Ponte05}, proven for general Lie quasi-bialgebroids.

We refer to $(M,\mathbf{a}_A\circ\rho,\pi_\DIF)$ as the \emph{induced quasi-Poisson $\g$-structure}.

\begin{remark}\label{QPBAisQPM}
As a converse to Example~\ref{QPMisQPBA}, suppose that $(TM,\rho,\DIF)$ is a quasi-Poisson $\g$-bialgebroid. Then by \cite[Lemma 2.2.]{GradedSymplSup}, $\DIF=[\pi,\cdot]$ for a unique bivector field $\pi\in\mathbf{\Gamma}(\wedge^2 TM)$. Since $\pi_\DIF^\sharp(df)=\mathbf{a}_A\circ \DIF f$, it follows that $\pi=\pi_\DIF$.

Consequently $(TM,\rho,\DIF)$ is of the form given in Example~\ref{QPMisQPBA} for the quasi-Poisson $\g$-structure $(M,\rho,\pi_\DIF)$.
\end{remark}

\begin{proposition}\label{liealgstruct2}
If $(A,\rho,\DIF)$ is a quasi-Poisson $\g$-bialgebroid, then $A^*$ becomes a Lie algebroid, where the Lie algebroid differential $d_{A^*}$ on $\mathbf{\Gamma}(\wedge A)$ is
\begin{equation}\label{LieAlgDiff2}d_{A^*}=\DIF+\frac{1}{2}\sum_i\rho(e^i)\wedge[\rho(e_i),\cdot]_A\,.\end{equation}
Furthermore, the action of $\g$ on $A^*$ preserves the Lie algebroid structure, and $\mu_\rho:A^*\to \g$ is a Lie algebroid morphism, where 
$\mu_\rho:=\on{pr}_\g\circ\rho^*$.
\end{proposition}

 A proof of this is given in \S~\ref{qpgbr}.

A quasi-Poisson $\g$-bialgebroid will be called \emph{integrable} if $A^*$ is an integrable Lie algebroid. In particular, $(A,\rho,\DIF)$ may be integrable even if $A$ is not.

We can now state our second theorem.

\begin{theorem}\label{mainthm2}
There is a one-to-one correspondence between source 1-connected Hamiltonian quasi-Poisson $\g$-groupoids, $(\Gamma,\rho_\Gamma,\pi_\Gamma,\Phi)$, and integrable quasi-Poisson $\g$-bialgebroids $(A,\rho,\DIF)$. Under this correspondence, the Lie algebroid of $\Gamma$ is $A^*$ and $\Phi$ integrates the Lie algebroid morphism $\mu_\rho:A^*\to\g$. 
\end{theorem}

We also have:

\begin{proposition}\label{maithm2pt2}
Suppose $(\Gamma,\rho_\Gamma,\pi_\Gamma,\Phi)$ is a Hamiltonian quasi-Poisson $\g$-groupoid corresponding to the quasi-Poisson $\g$-bialgebroids $(A,\rho,\DIF)$. Then the source map $s:\Gamma\to M$ is a quasi-Poisson morphism onto the induced quasi-Poisson $\g$-structure $(M,\mathbf{a}\circ\rho,\pi_\DIF)$ described in Proposition~\ref{inducedqpstruct}. Meanwhile the target map $t:\Gamma\to M$ is anti-quasi-Poisson.
\end{proposition}

We will provide a proof of both theorems in the next section using graded Poisson-Lie groups. 

\section{Quasi-Poisson structures and graded Poisson geometry}\label{part1}
In this section we make use of graded geometry (super geometry) to prove the results from \S~\ref{sect1}. Some good references for super geometry are \cite{DM99,Varadarajan04,Voronov91,VMP90,Kostant77}. For the additional structure of graded manifolds one may look at \cite{Severa01,Voronov02,Mehta06,GradedSymplSup}. 

%
%

\subsection{$\g$-differential algebras}\label{gdiffpre}
Let $\g$ be a Lie algebra, and $\hat\g=\g[1]\oplus\g\oplus \mathbb{R}[-1]$ be the graded Lie algebra with bracket given by
\begin{align*}
 [I_\xi,I_\eta]&=0  \\
 [L_\xi,I_\eta]&=I_{[\xi,\eta]_\g} &[L_\xi,L_\eta]&=L_{[\xi,\eta]_\g}  \\
 [D,I_\eta]&=L_\eta &[D,L_\eta]&=0 &[D,D]&=0
\end{align*}
Here $D$ is the generator of $\mathbb{R}[-1]$, and $L_\xi\in \g\subset\hat\g$ and $I_\xi\in \g[1]\subset\hat\g$ denote the elements corresponding to $\xi\in\g$.

If $\rho:\g\to\mathbf{\Gamma}(TM)$ is a morphism of Lie algebras, then the graded Lie algebra $\hat\g$ acts on the graded algebra $\Omega(M)$ by derivations.  $D$ acts by the de Rham differential $d$, $L_\xi$ acts by the Lie derivative $\Lied_{\rho(\xi)}$, and $I_\xi$ acts by the interior product $\iota_{\rho(\xi)}$.

Generally, a graded algebra with a graded action of $\hat\g$ by derivations is called a $\g$-differential algebra \cite{EquivCohom,GS99}. Note that, by a graded action, we mean that a degree $k$ element of $\hat\g$ acts by a degree $k$ derivation.

\begin{example}\label{gdiffalgex}
As a generalization of $\Omega(M)$, suppose $A\to M$ is any Lie algebroid, and $\rho:\g\to \mathbf{\Gamma}(A)$ is any Lie algebra morphism. Then $\mathbf{\Gamma}(\wedge A^*)$ is a $\g$-differential algebra. The action of $\hat\g$ is given as follows 
\begin{itemize}
\item $D\cdot \alpha=d_A\alpha$, where $\alpha\in\mathbf{\Gamma}(\wedge A^*)$ and $d_A$ is the Lie algebroid differential.
\item $I_\xi\cdot\alpha=\iota_{\rho(\xi)}\alpha$ for any $\xi\in\g$.
\item $L_\xi\cdot\alpha=\iota_{\rho(\xi)}(d_A\alpha)+d_A(\iota_{\rho(\xi)}\alpha)$ for $\xi\in\g$ and $\alpha\in\mathbf{\Gamma}(\wedge A^*)$.
\end{itemize}
\end{example}

\begin{remark}\label{supergeointerp}
If $A\to M$ is any vector bundle, the following are equivalent:
\begin{itemize}
\item $A\to M$ is a Lie algebroid, and there is a Lie algebra morphism $\rho:\g\to \mathbf{\Gamma}(A)$
\item $\mathbf{\Gamma}(\wedge A^*)$ is a $\g$-differential algebra.
\end{itemize}

We can think of $\mathbf{\Gamma}(\wedge A^*)$ as the algebra of functions on the graded manifold $A[1]$, hence another equivalent formulation is
\begin{itemize}
\item The graded Lie algebra $\hat\g$ acts on the graded manifold $A[1]$.
\end{itemize}

\end{remark}

\subsection{The quadratic graded Lie algebra $\mathcal{Q}(\g)$}\label{QGdef}
Suppose a Lie algebra $\g$  possesses an invariant non-degenerate symmetric bilinear form $\langle\cdot,\cdot\rangle_\g$. We can associate to $\g$ the quadratic graded Lie algebra $\mathcal{Q}(\g)$, an $\mathbb{R}[2]$ central extension  of $\hat\g$, $$0\to\mathbb{R}[2]\to\mathcal{Q}(\g)\to\hat\g\to 0,$$  which plays a central role in the theory of quasi-Poisson structures.

As a graded vector space, $\mathcal{Q}(\g)=\mathbb{R}[2]\oplus\hat\g$. Let $T$ denote the generator of $\mathbb{R}[2]$. The central extension is given by the cocycle 
$$c(I_u,I_v)=\langle u,v\rangle T,\quad c(D,\cdot)=c(L_u,\cdot)=0,$$ i.e.\ $$[a,b]_{\mathcal{Q}(g)}=[a,b]_{\hat\g}+c(a,b)$$
for $a,b\in\hat\g$.

The quadratic form $\langle\cdot,\cdot\rangle_{\mathcal{Q}(\g)}$ of degree $1$ is given by 
\begin{subequations}\begin{equation}\label{deg1cond}\langle a,b\rangle_{\mathcal{Q}(\g)}=0\text{ for any } a,b\in\mathcal{Q}(\g)\text{ such that }\on{deg}(a)+\on{deg}(b)+1\neq0,\text{ and}\end{equation}
\begin{equation}\langle T,D\rangle_{\mathcal{Q}(\g)}=1,\quad \langle I_\xi,L_\eta\rangle_{\mathcal{Q}(\g)}=\langle\xi,\eta\rangle_\g.\end{equation}\end{subequations} 
Note that \eqref{deg1cond} is equivalent to saying the quadratic form is of degree $1$.

\begin{remark}The Lie algebra $\mathcal{Q}(\g)$ was first introduced in \cite{Alekseev05}, where the so called non-commutative Weil algebra was defined as a quotient of the enveloping algebra of $\mathcal{Q}(\g)$.
\end{remark}


\subsection{Quasi-Poisson $\g$-manifolds revisited}\label{mathbbR2timesg1}
It is easy to check that $\mathbb{R}[2]\oplus\g[1]$ and $\g\oplus\mathbb{R}[-1]$ are transverse Lagrangian subalgebras of $\mathcal{Q}(\g)$. 
Therefore $(\mathcal{Q}(\g),\mathbb{R}[2]\oplus\g[1],\g\oplus\mathbb{R}[-1])$ forms a Manin triple \cite{Drinfeld83,thesis-3,lu90}. The corresponding Lie bialgebra $\mathbb{R}[2]\oplus\g[1]$ integrates to the Poisson Lie group $$\mathbf{G}_{small}=\mathbb{R}[2]\times\g[1],$$ where multiplication is given by
$$(t,\xi)\cdot(t',\xi')=(t+t'+\frac{1}{2}\langle\xi,\xi'\rangle_\g,\xi+\xi').$$

 Since the quadratic form on $\mathcal{Q}(\g)$ is of degree $1$, the Poisson bracket on $\mathbf{G}_{small}$ is of degree $-1$. To describe the Poisson bracket, note that linear functions on $\g[1]$ may be identified with elements of $\g$ (using the quadratic form). If we let $t$ denote the standard coordinate on $\mathbb{R}[2]$ then we see that there is a canonical algebra isomorphism $C^\infty(\mathbb{R}[2]\times\g[1])\cong (\wedge^*\g)[t]$. Under this isomorphism the Poisson bracket is simply
$$\{t,t\}=\phi \quad \{t,\xi\}=0 \quad \{\xi,\eta\}=[\xi,\eta]_\g$$
\begin{proposition}\label{prop:alt_qpoiss}
A quasi-Poisson $\g$-structure on $M$ is equivalent to a graded Poisson map $T^*[1]M\to\mathbf{G}_{small}$.
\end{proposition}
\begin{proof}
A quasi-Poisson $\g$-manifold $(M,\rho,\pi)$ is equivalent to a morphism of Gerstenhaber algebras $\rho':(\wedge^*\g)[t]\to \mathbf{\Gamma}(\wedge^* TM)$. Here $\rho'$ is defined on the generators by $\rho'(t)=\pi$ and $\rho'(\xi)=\rho(\xi)$ for any $\xi\in\g$. 

The standard symplectic form on $T^*[1]M$ induces a degree $-1$ Poisson bracket on $C^\infty(T^*[1]M)$. As Gerstenhaber algebras $C^\infty(T^*[1]M)\cong\mathbf{\Gamma}(\wedge^*TM),$ canonically. Paired with the isomorphism $C^\infty(\mathbf{G}_{small})\cong (\wedge^*\g)[t]$, we see that $\rho'$ defines a morphism of Poisson algebras $$C^\infty(\mathbf{G}_{small})\to C^\infty(T^*[1]M).$$ This is equivalent to a Poisson morphism $$T^*[1]M\to \mathbf{G}_{small}.$$
\end{proof}

\begin{proof}[Proof of Theorem~\ref{liealgstruct}]
If $(M,\rho,\pi)$ is a quasi-Poisson $\g$-manifold, then we have a Poisson map $f:T^*[1]M\to\mathbf{G}_{small}$ ($f^*t=\pi$, $f^*\xi=\rho(\xi)$). Therefore, the dual Lie algebra $\g\oplus\mathbb{R}[-1]$ acts on $T^*[1]M$. To describe the action explicitly, recall \cite{thesis-3} that the left invariant one forms on $\mathbf{G}_{small}$ form a subalgebra of $\mathbf{\Gamma}(T^*\mathbf{G}_{small})$ isomorphic to $\g\oplus\mathbb{R}[-1]\cong T^*_e[1]\mathbf{G}_{small}$ (evaluation at the identity provides the isomorphism). The left-invariant 1-form on $\mathbf{G}_{small}$ corresponding to $D$ is 
$$dt+\frac{1}{2}\sum_i\xi^id\xi_i,$$
where $\xi_i$ and $\xi^i$ refer to coordinates on $\g[1]$ induced by the basis vectors $e_i$ and $e^i$, respectively.
The corresponding vector field on $T^*[1]M$ is thus
$$\{f^*t,\cdot\}+\frac{1}{2}\sum_if^*\xi^i\{f^*\xi_i,\cdot\},$$
i.e.\ the differential $d_{T^*M}$ \eqref{liealgdiff}. Since $[D,D]=0$, this shows that $d_{T^*M}^2=0$. The action of $\g$ on $T^*[1]M$ preserves $d_{T^*M}$ (since $\g\oplus\mathbb{R}[-1]$ is a direct sum) and it is just the natural lift of the action $\rho$ on $M$ (the left-invariant 1-form on $\mathbf{G}_{small}$ corresponding to $\xi\in\g$ is $d\xi$).

The dressing action of $D\in\hat\g$ on $\mathbf{G}_{small}$ is given by
$$\{t,\cdot\}+\frac{1}{2}\sum_i\xi^i\{\xi_i,\cdot\}=\phi\partial_t+d_\g,$$
where $d_\g$ is the Lie algebra differential of $\g$. The projection $\mathbf{G}_{small}\to\g[1]$ is thus $\mathbb{R}[-1]$-equivariant, with respect to the $\mathbb{R}[-1]$ actions generated by $D\in\hat\g$ and $d_\g$, respectively. Since the map $f$ is also $\mathbb{R}[-1]$-equivariant, so is their composition $T^*[1]M\to\g[1]$, i.e.\ we have a Lie algebroid morphism $T^*M\to\g$.
\end{proof}

The fusion also appears in a natural way from this perspective. A Poisson morphism $$f:T^*[1]M\to \mathbf{G}_{small}\times\mathbf{G}_{small}$$ defines a quasi-Poisson $\g\oplus\g$-structure on $M$. Since $\mathbf{G}_{small}$ is a Poisson Lie group, the multiplication map $$\on{mult}:\mathbf{G}_{small}\times\mathbf{G}_{small}\to\mathbf{G}_{small}$$ is a Poisson morphism. The map $$\on{mult}\circ f:T^*[1]M\to \mathbf{G}_{small}$$ defines a quasi-Poisson $\g$-structure on $M$. Since $\on{mult}^*t=t_1+t_2+\frac{1}{2}\sum_i(\xi^i)_1(\xi_i)_2$ (where the sub-indices $(\cdot)_1$ and $(\cdot)_2$ indicate which factor of $\mathbf{G}_{small}\times\mathbf{G}_{small}$ the coordinates parametrize), the bivector on $M$ is modified by the term $f^*(\frac{1}{2}\sum_i(\xi^i)_1(\xi_i)_2)$ (note the similarity to \eqref{psi}). It is easy to check that this is the same quasi-Poisson $\g$-structure on $M$ given by the fusion \eqref{fusion}.

\subsection{Hamiltonian quasi-Poisson $\g$-manifolds revisited}\label{qhqpgmr}
Let $G$ be a Lie group with Lie algebra $\g$. Suppose $\langle\cdot,\cdot\rangle_\g$ is a quadratic form on $\g$. Let $\bar\g$ denote the quadratic Lie algebra whose quadratic form is $\langle\cdot,\cdot\rangle_{\bar\g}=-\langle\cdot,\cdot\rangle_\g$.
We let $\dd=\g\oplus\bar\g$, and $\on{diag}(\g)\subset\dd$ denote the  diagonal subalgebra. 
\begin{subequations}\label{eq:ghat_bialg}
Then 
\begin{equation}\label{eq:ghat_bialga}
\mathbb{R}[2]\oplus\g[1]\oplus\bar\g
\end{equation}
 and 
 \begin{equation}\label{eq:ghat_bialgb}
\hat\g\cong\on{diag}(\g)[1]\oplus\on{diag}(\g)\oplus\mathbb{R}[-1]
\end{equation}
\end{subequations}
are two Lagrangian subalgebras of the quadratic Lie algebra $$\mathcal{Q}(\dd)=\mathbb{R}[2]\times\dd[1]\times\dd\times\mathbb{R}[-1]$$ defined in \S~\ref{QGdef}. The corresponding Lie bialgebra $\mathbb{R}[2]\oplus\g[1]\oplus\bar\g$ integrates to the Poisson Lie group $$\mathbf{G}_{big}=\mathbb{R}[2]\times\g[1]\times G,$$ where multiplication is given by
$$(t,\xi,g)\cdot(t',\xi',g')=(t+t'+\frac{1}{2}\langle\xi,\xi'\rangle,\xi+\xi',g\cdot g').$$
 The group $\mathbf{G}_{big}$ is thus the direct product of $G$ with the Heisenberg group $\mathbf{G}_{small}$ described above. As in \S~\ref{mathbbR2timesg1}, there is a canonical identification $C^\infty(\mathbf{G}_{big})\cong(\wedge^*\g)[t]\otimes C^\infty(G)$. Using this identification, we may describe the Poisson bracket (of degree $-1$) on $\mathbf{G}_{big}$ by
\begin{align}
%
\{t,t\}&=\phi\\
\{t,\xi\}&=0 & \{\xi,\eta\}&=[\xi,\eta]_\g\\
\label{bfGpstruct} \{t,f\}&=b^*df & \{\xi,f\}&=(\xi^L-\xi^R)\cdot f & \{f,g\}&=0 
\end{align}
where $f,g\in C^\infty(G)$, $\xi,\eta\in\g$, $\xi^L$ and $\xi^R$ denote the corresponding left and right invariant vector fields on $G$ and $b$ is given by \eqref{beq}.


We have
\begin{proposition}\label{bfGmap}
A Hamiltonian quasi-Poisson $\g$-structure on $M$ is equivalent to a graded Poisson map $T^*[1]M\to\mathbf{G}_{big}$.
\end{proposition} 

\begin{proof}
The proof of Proposition~\ref{prop:alt_qpoiss} shows that $M$ is a quasi-Poisson $\g$-manifold. The map $T^*[1]M\to\mathbf{G}_{big}$ restricts to define a map \begin{equation}\label{restrmommap}\Phi:M\to G;\end{equation} and the formulas for the brackets in \eqref{bfGpstruct} show that $\Phi$ defines a moment map (Definition~\ref{MomentMapDefn}).
\end{proof}

As in \S~\ref{mathbbR2timesg1}, fusion can be described in terms of composition with the multiplication Poisson morphism $$\on{mult}:\mathbf{G}_{big}\times\mathbf{G}_{big}\to\mathbf{G}_{big},$$ this is precisely equivalent to the explanation given in \cite{PureSpinorsOnL}.

\begin{proof}[Proof of Theorem~\ref{g-inclusion}]
A Poisson morphism $F:T^*[1]M\to\mathbf{G}_{big}$ induces an action of $\on{diag}(\g)[1]\oplus\on{diag}(\g)\oplus\mathbb{R}[-1]\cong\hat\g$ on $T^*[1]M$, $$\varrho:\hat\g\to\mathbf{\Gamma}\big(T(T^*[1]M)\big).$$ For $\xi\in\g$, and the corresponding element $I_\xi\in\hat\g$, let us describe this action explicitly.

 The standard symplectic form on $T^*[1]M$ induces a degree $-1$ Poisson bracket on $C^\infty(T^*[1]M)$. As Gerstenhaber algebras $C^\infty(T^*[1]M)\cong\mathbf{\Gamma}(\wedge^*TM),$ canonically. A function $f\in C^\infty(M)$ acts on $\mathbf{\Gamma}(\wedge^*TM)$ by contraction with $df$. It follows that the Hamiltonian vector field generated by a one form $\alpha\in\Omega^1(M)$ acts on $\mathbf{\Gamma}(\wedge^*TM)$ by contraction with $\alpha$.

Let $\lambda_{I_\xi}$ denote the left invariant one form on $\mathbf{G}_{big}$ corresponding to $I_\xi$. $\lambda_{I_\xi}$ is just the pullback of $\langle\xi,\theta^L\rangle\in\Omega^1(G)$ to $\mathbf{G}_{big}$. Consequently, $\varrho(I_\xi)$ acts by contraction with $\Phi^*\langle\xi,\theta^L\rangle$, where $\Phi:M\to G$ given by \eqref{restrmommap}. For $X\in\mathbf{\Gamma}(TM)$,
$$\varrho(I_\xi)\cdot X=\bigl(\Phi^*\langle\xi,\theta^L\rangle\bigr)(X)=i(\xi)(X),$$
where $i:\g\to \Omega^1(M)$ is given by \eqref{i-map-def}.

As stated in Remark~\ref{supergeointerp}, an action of $\hat\g$ on $T^*[1]M$ is equivalent to $\mathbf{\Gamma}(\wedge^* TM)$ being a $\g$-differential algebra. The proof of Theorem~\ref{liealgstruct} shows that the differential is given by \eqref{liealgdiff}. From this perspective, Theorem~\ref{g-inclusion} is just a special case of Example~\ref{gdiffalgex}.
\end{proof}

\subsection{Hamiltonian quasi-Poisson $\g$-groupoids revisited}\label{qhqpgr}Let $\Gamma$ be any groupoid. Recall that $T^*[1]\Gamma$ has a natural groupoid structure (see Appendix~A, Page~\pageref{vbgroupoids}). Combining this structure with the canonical symplectic structure on the cotangent bundle, $T^*[1]\Gamma$ becomes a symplectic groupoid.

Since $T^*[1]\Gamma$ is a symplectic groupoid, it integrates a Poisson manifold \cite{weinstein87-1}. We can describe this Poisson manifold explicitly. Let $A^*$ denote Lie algebroid corresponding to the groupoid $\Gamma$ (we denote the Lie algebroid $A^*$ (and not $A$) for later convenience). $A[1]$ has a linear Poisson structure on it (of degree $-1$) defining the Lie algebroid structure on $A^*$ \cite{LieAlgebroidsH}. $T^*[1]\Gamma$ is the symplectic groupoid integrating $A[1]$.

\begin{proposition}\label{prop:qhqpgr}
A compatible Hamiltonian quasi-Poisson $\g$-structure on $\Gamma$ is equivalent to a morphism of Poisson groupoids 

\begin{equation}\label{qhqpgroupoid}F:T^*[1]\Gamma\to\mathbf{G}_{big}.\end{equation}
\end{proposition}
\begin{proof}
First we introduce some notation. If $M$ is a graded manifold and $f\in C^\infty(M)$, let $(f)_i$ denote the pullback of $f$ to the $i^{th}$ factor of the direct power $M^n$. If $M$ is a graded Poisson manifold with Poisson bracket $\{\cdot,\cdot\}_M$, let $\bar M$ denote the same graded manifold with the Poisson bracket \begin{equation}\label{barM}\{\cdot,\cdot\}_{\bar M}=-\{\cdot,\cdot\}_M.\end{equation}

By Proposition~\ref{bfGmap}, $F$ defines a Hamiltonian quasi-Poisson $\g$-structure on $\Gamma$. The moment map $\Phi:\Gamma\to G$ is given by restricting $F$ to the subgroupoid $\Gamma\subset T^*[1]\Gamma$. Consequently $\Phi$ is a morphism of Lie groupoids.

Under the isomorphism $C^\infty(\mathbf{G}_{big})\cong (\wedge^*\g)[t]\otimes C^\infty(G)$, $\eta\in\g$ defines a function on $\mathbf{G}_{big}$. We notice that the functions $$(\eta)_1+(\eta)_2-(\eta)_3,\;\text{ and }\; \big(t_1+t_2+\frac{1}{2}(\xi^i)_1(\xi_i)_2\big)-t_3$$ vanish on the graph of the multiplication $\on{Gr}_{\on{mult}_{\mathbf{G}_{big}}}\in\mathbf{G}_{big}\times\mathbf{G}_{big}\times\overline{\mathbf{G}_{big}}$.  Since $F$ is a groupoid morphism, it follows that the functions
$$(F^*\eta)_1+(F^*\eta)_2-(F^*\eta)_3,\;\text{ and }\; \big((F^*t)_1+(F^*t)_2+\frac{1}{2}\sum_i(F^*\xi^i)_1(F^*\xi_i)_2\big)-(F^*t)_3$$
vanish on the graph of the multiplication $\on{Gr}_{\on{mult}_{T^*[1]\Gamma}}$. In the first case, this shows that action of $\g$ on $\Gamma\times\Gamma\times\Gamma$ is tangent to the graph of the multiplication. In other words, $\g$ acts on $\Gamma$ by groupoid automorphisms. In the latter case, this shows that the bivector field on $\Gamma$ is fusion multiplicative.

\end{proof}

\subsection{The Lie bialgebra $\hat{\mathfrak{g}}$ is quasi-triangular}\label{sect:ghat-qtr}
Recall that \eqref{eq:ghat_bialg} makes $\hat\g$ into a Lie bialgebra.

\begin{proposition}\label{ghat_is_qtriang}
The element $\hat r=\sum_iI_{e^i}\otimes L_{e_i}\in\hat\g\otimes\hat\g$ is an $r$-matrix for the graded Lie bialgebra $\hat\g$.
\end{proposition}
\begin{proof}
We may view the degree $-1$ element $\hat r\in\hat\g\otimes\hat\g$ as a degree $0$ map $\hat r:\hat\g^*[1]\to\hat\g$.
Using \eqref{eq:ghat_bialg}, we identify $\hat\g^*[1]$ and $\hat\g$ with the transverse Lagrangian subalgebras  $$\mathbb{R}[2]\oplus\g[1]\oplus\bar\g\subset\mathcal{Q}(\dd)$$ and $$\on{diag}(\g)[1]\oplus\on{diag}(\g)\oplus\mathbb{R}[-1]\subset\mathcal{Q}(\dd)$$ respectively.  Then the graph of $\hat r:\hat\g^*[1]\to\hat\g$ is identified with the subspace $$\on{Gr}(\hat r)\cong\mathbb{R}[2]\oplus\g[1]\oplus\g.$$ Thus $\on{Gr}(\hat r)$ is an ideal of the Drinfeld double of $\hat\g$. Equivalently $\hat r$ is an $r$-matrix.
%
%
%
%
%
\end{proof}

\begin{proof}[Proof of Proposition \ref{prop:lalg-mon}]

By Proposition~\ref{ghat_is_qtriang}, $\hat\g$ is a quasi-triangular Lie bialgebra. Recall from \S~\ref{FusionSect}, the functor $$F_{\hat\g}:\on{Ham-qPois}_{\hat\g}\to\on{Ham-Pois}_{\hat\g},$$ from the category of Hamiltonian quasi-Poisson $\hat\g$-manifolds to the category of Hamiltonian Poisson $\hat\g$-manifolds, and its inverse $F^{-1}_{\hat\g}$, which describe an equivalence of categories.

Let $M$ and $M'$ be two Hamiltonian quasi-Poisson $\g$-manifolds. By Proposition~\ref{bfGmap}, $T^*[1]M$ and $T^*[1]M'$ are Hamiltonian Poisson $\hat\g$-manifolds. Now $F_{\hat\g}$ is strongly monoidal, and we have the natural transformation
\begin{equation}\label{eq:nattranshatg}\mathcal{J}:F_{\hat\g}\big(F^{-1}_{\hat\g}(T^*[1]M)\circledast F^{-1}_{\hat\g}(T^*[1]M')\big)\to (T^*[1]M)\times(T^*[1]M')\end{equation}
 described in Remark~\ref{rem:mmequiv}. Note that the map $j:\mathbf{G}_{big}\to\hat G$ from Remark~\ref{rem:mmequiv} is just the projection $\mathbf{G}_{big}\to\g[1]$. Therefore \eqref{eq:nattranshatg} is just the map $J$ described in Proposition~\ref{prop:lalg-mon}.

The right hand side of \eqref{eq:nattranshatg} describes the fusion of $T^*[1]M$ and $T^*[1]M'$ as Hamiltonian Poisson $\hat\g$-manifolds. As shown in \S~\ref{qhqpgmr}, $(T^*[1]M)\times(T^*[1]M')=T^*[1](M\circledast M')$, where the right hand side is the fusion of $M$ and $M'$ as Hamiltonian quasi-Poisson $\g$-manifolds.

The left hand side of $\eqref{eq:nattranshatg}$ describes the fusion of $F^{-1}_{\hat\g}(T^*[1]M)$ and $F^{-1}_{\hat\g}(T^*[1]M')$ as Hamiltonian quasi-Poisson $\hat\g$-manifolds. Therefore the action of $\hat\g$ on $(T^*[1]M)\circledast(T^*[1]M')$ is by definition diagonal, and the functor $F_{\hat\g}$ preserves this action. 

Therefore the Lie algebroid structure on the left-hand side (given by the action of $D\in\hat\g$) is the direct sum of the Lie algebroids $T^*M$ and $T^*M'$ (as the action of $\hat\g$ is diagonal), while the right-hand side corresponds to the Lie algebroid $T^*(M\circledast M')$.




\end{proof}

\subsection{Quasi-Poisson $\g$-bialgebroids revisited}\label{qpgbr}
Before proving Theorem~\ref{mainthm2}, it is important to formulate a description of quasi-Poisson $\g$-bialgebroids in terms of the Manin triple $(\mathcal{Q}(\g\oplus\bar\g),\hat\g,\mathbb{R}[2]\oplus\g[1]\oplus\bar\g)$. In fact, the description is quite natural, namely:

\begin{proposition}\label{p-qpgb-htpa}
Suppose $A\to M$ is a vector bundle.
The following are equivalent:
\begin{itemize}
\item $A$ is a quasi-Poisson $\g$-bialgebroid.
\item There is a Poisson structure of degree $-1$ on $A[1]$, and a Lie bialgebra action of $\hat\g$ on $A[1]$.
\end{itemize}
\end{proposition}

\begin{proof}

Suppose we are given a degree $-1$ Poisson structure $\pi$ on $A[1]$ and an action $\hat\rho$ of the Lie bialgebra $\hat\g$ on $A[1]$. We need to show that $A$ is a quasi-Poisson $\g$-bialgebroid. Recalling Definition~\ref{def:qpgb}, we must show that

\begin{itemize}
\item {\em We have a Lie algebroid structure on $A$ and a Lie algebra morphism $$\rho:\g\to\mathbf{\Gamma}(A).$$}
\end{itemize}

(This follows directly from Remark~\ref{supergeointerp}.)

\begin{itemize}
\item{\em There is a degree $+1$ derivation $\DIF$ of the Gerstenhaber algebra $\mathbf{\Gamma}(\bigwedge A)$, such that 
\begin{itemize}
\item $\DIF\rho(\xi)=0$ for any $\xi\in\g$, and
\item $\DIF^2=\frac{1}{2}[\rho(\phi),\cdot]\,$.
\end{itemize}}
\end{itemize}

By Proposition~\ref{ghat_is_qtriang} and the graded version of Remark~\ref{qaction_vs_bialg}, an action of the Lie bialgebra $\hat\g$ on $A[1]$ is equivalent to a quasi-Poisson action of $\hat\g$ on $A[1]$. Let us describe this explicitly, to avoid possible sign problems (as $\hat\g$ is a graded Lie bialgebra with cobracket of degree $-1$). A Poisson structure of degree $-1$ on the graded manifold $A[1]$ is, by definition, a function $\pi$ on the bigraded symplectic manifold $T^*[1,1]A[1,0]$ of degree $(1,2)$ such that $\{\pi,\pi\}=0$. An action $\hat\rho$ of the graded Lie algebra $\hat\g$ can be seen as a map $\hat\rho:\hat\g\to C^\infty(T^*[1,1]A[1,0])$ (as vector fields can be seen as linear functions on the cotangent bundle) shifting degrees by (1,1). 

The action $\hat\rho$ is a Lie bialgebra action on $(A[1,0],\pi)$. Therefore, by Proposition~\ref{ghat_is_qtriang} and Remark~\ref{qaction_vs_bialg},
\begin{equation}\label{LieAlgDiffEq}\tilde\pi=\pi-\frac{1}{2}\sum_i\hat\rho(I_{e^i})\hat\rho(L_{e_i})\end{equation}
is $\hat\g$-invariant and $(A[1,0],\tilde\pi)$ is a quasi-Poisson $\hat\g$-space:
$$\{\tilde\pi,\tilde\pi\}=\frac{1}{4}\sum_{ijk}c_{ijk}\, \hat\rho(I_{e^i})\hat\rho(I_{e^j})\hat\rho(L_{e^k}),$$
where  $c_{ijk}=\langle[e_i,e_j],e_k\rangle$ are the structure constants of $\g$. We can rewrite it as
\begin{equation}\label{eq:pitilde2}
\{\tilde\pi,\tilde\pi\}=\{\hat\rho(D),\hat\rho(I_\phi)\}.
\end{equation}

Using the canonical symplectomorphism
\begin{equation}\label{symplmorph}T^*[1,1]A[1,0]\cong T^*[1,1]A^*[0,1],\end{equation}
$\tilde\pi$ becomes a {\em vector field} on $A^*[0,1]$ (since it is a function linear on the fibers of $T^*[1,1]A^*[0,1]$), i.e.\ a derivation $\DIF$ of the algebra $\mathbf{\Gamma}(\wedge A)$ of degree 1. Since $\tilde\pi$ is $\hat\g$-invariant, $\DIF$ preserves the Gerstenhaber bracket on $\mathbf{\Gamma}(\wedge A)$ and $\DIF\rho(\xi)=0$ for every $\xi\in\g$. Finally, Equation \eqref{eq:pitilde2} becomes $\DIF^2=\frac{1}{2}[\rho(\phi),\cdot]$.

 We have shown that, $(A,\DIF,\rho)$ is a quasi-Poisson $\g$-bialgebroid. To establish the converse, just reverse the procedure.

\end{proof}

\begin{proof}[Proof of Proposition~\ref{liealgstruct2}]
 By Proposition~\ref{p-qpgb-htpa}, a quasi-Poisson $\g$-bialgebroid structure on $A$ defines a degree $-1$ Poisson structure on $A[1]$. This is equivalent to a Lie algebroid structure on $A^*$ \cite{LieAlgebroidsH}.

 A careful examination of the proof of Proposition~\ref{p-qpgb-htpa} allows us to describe the Lie algebroid differential explicitly. The bivector corresponding to the Poisson structure on $A[1]$ is a function $\pi$ on the bigraded symplectic manifold $T^*[1,1]A[1,0]$ of degree $(1,2)$. Under the canonical symplectomorphism \eqref{symplmorph} it becomes a degree $+1$ vector field on $A^*[1]$.  By \eqref{LieAlgDiffEq} this vector field is $$\DIF+\frac{1}{2}\sum_i\hat\rho(I_{e^i})\hat\rho(L_{e_i}).$$ The corresponding Lie algebroid differential is given by \eqref{LieAlgDiff2}.

The remaining details in Proposition~\ref{liealgstruct2} follow from a similar examination of the proof of Proposition~\ref{p-qpgb-htpa}.
\end{proof}

\begin{remark}\label{rem:alt_qPbialg}
One can also describe quasi-Poisson $\g$-bialgebroids in the spirit of Proposition~\ref{prop:alt_qpoiss}. If $A$ is a Lie algebroid, so that $A^*[1]$ is a graded Poisson manifold, a Poisson map $A^*[1]\to\mathbf{G}_{small}$ would give us a quasi-Poisson $\g$-bialgebroid structure on $A$, but with the additional property that $\DIF$ is Hamiltonian. In general, a quasi-Poisson $\g$-bialgebroid structure on a Lie algebroid $A$ is equivalent to a principal Poisson $\mathbb{R}[2]$-bundle $P\to A^*[1]$ with a Poisson $\mathbb{R}[2]$-equivariant map $P\to\mathbf{G}_{small}$.
\end{remark}

\subsection{Proof of Theorem~\ref{mainthm2}}\label{proofthm1}

The proof of \cite[Theorem 5.5]{Xu95} (see also \cite{Fernandes07,Fernandes09}) goes through in the graded setting to show that the existence of a morphism of Poisson groupoids
$$F:T^*[1]\Gamma\to\mathbf{G}_{big}$$
is equivalent to the action of the Lie bialgebra $\hat\g$ on the Poisson manifold $A[1]$. By Proposition \ref{prop:qhqpgr}, the former  describes a compatible Hamiltonian quasi-Poisson $\g$-structure on $\Gamma$, while the latter describes a quasi-Poisson $\g$-bialgebroid structure $(A,\rho,\DIF)$  (see Proposition~\ref{p-qpgb-htpa}). This proves Theorem~\ref{mainthm2}.

\begin{proof}[Proof of Proposition~\ref{maithm2pt2}]
Let $(\Gamma,\rho_\Gamma,\pi_\Gamma)$ be the quasi-Poisson structure on $\Gamma$, and $(M,\rho_M,\pi_\DIF)$ be the quasi-Poisson structure on $M$ induced by $(A,\rho,\DIF)$ (see Proposition~\ref{inducedqpstruct}). We must show that the source map $s_0:\Gamma\to M$ is a quasi-Poisson morphism. Let $D_A$ and $D_\Gamma$ be the homological vector fields on $A[1]$ and $T^*[1]\Gamma$ defined by the respective actions of $D\in\hat\g$. Since the source map $s:T^*[1]\Gamma\to A[1]$ is $\hat\g$-equivariant, we have
$$s^*(D_A f)=D_\Gamma s^*f,$$ for every $f\in C^\infty(M)$, where we view $f$ as an element of $C^\infty(A[1])$. Since $s:T^*[1]\Gamma\to A[1]$ is also a Poisson map, it follows that 
\begin{equation}\label{inveq}s^*\{D_A f,g\}_{A[1]}=\{D_\Gamma s^*f,s^*g\}_{T^*[1]\Gamma},\end{equation}
for every $g\in C^\infty(M)$. Now $D_A$ and $D_\Gamma$ define the Lie algebroid differentials on $\mathbf{\Gamma}(\wedge^* A^*)$ and $\mathbf{\Gamma}(\wedge^* T\Gamma)$ respectively, while $\{\cdot,\cdot\}_{A[1]}$ and $\{\cdot,\cdot\}_{T^*[1]\Gamma}$ define the respective Gerstenhaber algebra brackets, so the left hand side of \eqref{inveq} is equal to $$s^*\big[ \DIF f+\frac{1}{2}\sum_i\rho(e^i)\wedge[\rho(e_i),f],g\big],$$ or simply $$s_0^*\big[ [\pi_\DIF,f]+\frac{1}{2}\sum_i\rho_M(e^i)\wedge[\rho_M(e_i),f],g\big],$$ while the right hand side is $$\big[ [\pi_\Gamma,s_0^*f]+\frac{1}{2}\sum_i\rho_\Gamma(e^i)\wedge[\rho_\Gamma(e_i),s_0^*f],s_0^*g\big].$$
Since $s_0$ is $\g$-invariant, it follows that
$$s_0^*\big[ [\pi_\DIF,f],g\big]=\big[ [\pi_\Gamma,s_0^*f],s_0^*g\big].$$
Therefore, $s_0:(\Gamma,\rho_\Gamma,\pi_\Gamma)\to(M,\rho_M,\pi_\DIF)$ is a quasi-Poisson morphism. Similarily, one may check that the target map is anti-quasi-Poisson.
\end{proof}

\subsection{Proof of Theorem~\ref{mainthm1}}
Let $\Gamma\rightrightarrows M$ be a  Hamiltonian quasi-Poisson $\g$-groupoid.
By Proposition \ref{prop:lalg-mon} the Lie algebroid $T^*\Gamma$ is multiplicative, where the groupoid structure on $T^*\Gamma$ is the standard
one (see Appendix \ref{vbgroupoids}) precomposed with the map $J$. In particular, the anchor map $\mathbf{a}_\Gamma:T^*\Gamma\to T\Gamma$ is a groupoid morphism.

Suppose that $\Gamma$ is source 1-connected. Let $A\to M$ be the quasi-Poisson $\g$-bialgebroid corresponding to $\Gamma$. Then for any $x\in M\subset\Gamma$ we have $T_x\Gamma\cong T_xM\oplus A_x^*$ and $\mathbf{a}_\Gamma:A_x\oplus T_x{}^*M\to T_xM\oplus A_x^*$ equal to $\mathbf{a}_A\oplus \mathbf{a}^*_A$. Hence  $\mathbf{a}_A$ is an isomorphism if and only if the anchor map $\mathbf{a}_\Gamma$ is an isomorphism at points of $M$. However if $\mathbf{a}_\Gamma$ is an isomorphism at points of $M$, then it must an isomorphism in a neighborhood of $M$. Since $\mathbf{a}_\Gamma$ it is a morphism of source 1-connected groupoids, this is the case if and only if it is an isomorphism everywhere.

By Theorem~\ref{thm:crit_for_qH}, $\Gamma$ is a Hamiltonian quasi-symplectic $\g$-space if and only if $\mathbf{a}_\Gamma$ is an isomorphism. By Remark \ref{QPBAisQPM}, $A$ comes from a quasi-Poisson structure on $M$ if and only if $\mathbf{a}_A$ is an isomorphism. As we just proved, these two conditions are equivalent, i.e.\ Theorem~\ref{mainthm1} is proven.

\section{Examples}\label{examples1}

\subsection{Quasi-symplectic case}
Let $(M,\rho,\pi)$ be a quasi-Poisson $\g$-manifold such that its anchor $\mathbf{a}:T^*M\to TM$ is bijective, i.e.\ $M$ is a quasi-symplectic manifold. Since $\mathbf{a}$ is an isomorphism of Lie algebroids, the source-1-connected groupoid integrating $T^*M$ is the fundamental groupoid $\Pi(M)$ of $M$. 
$\Pi(M)$ is a covering of $M\times M$, and the quasi-Poisson structure on $\Pi(M)$ is the lift of the quasi-Poisson structure on $(M,\pi,\rho)\circledast(M,-\pi,\rho)$.

The Lie algebroid morphism $\mu_\rho:T^*M\to\g$ gives us (via $\mathbf{a}$) a Lie algebroid morphism $TM\to\g$, i.e.\ a flat $\g$-connection on $M$. The moment map $\Pi(M)\to G$ is the parallel transport of this connection.

Suppose $(M,\rho,\phi)$ is endowed with a moment map $\Phi$ (so that it is a Hamiltonian quasi-symplectic $\g$-manifold). In this case, a Hamiltonian quasi-symplectic groupoid integrating it is just the pair groupoid $(M,\rho,\pi,\Phi)\circledast(M,\rho,-\pi,\Phi^{-1})$; to get a source-1-connected groupoid, we just lift the Hamiltonian quasi-symplectic structure to $\Pi(M)$.

\begin{remark}[Poisson Parallel]
This example is the quasi-Poisson analogue of an example in Poisson geometry. Namely, if $(N,\pi_N)$ is a Poisson manifold so that the anchor $\pi_N^\sharp:T^*N\to TN$ is bijective, then a symplectic groupoid integrating $(N,\pi_N)$ is just the pair groupoid $(N,\pi)\times(N,-\pi)$. If $\Pi(N)$ is the fundamental groupoid of $N$, then $\Pi(N)$ is a covering of $N\times N$. Consequently it inherits the structure of a symplectic groupoid. As such $\Pi(N)$ is the source 1-connected symplectic groupoid integrating $N$.
\end{remark}

\subsection{The double}

For a related example, let $G$ be a Lie group with Lie algebra $\g$. For $\xi\in\g$, let $\xi^L$ and $\xi^R$ denote the corresponding left and right invariant vector fields on $G$, and let $\rho:\g\oplus\g\to \mathbf{\Gamma}(TG)$ be given by $\rho(\xi,\eta)=-\xi^R+\eta^L$. As explained in \cite[Example~5.3]{Alekseev00}, since $\rho(\phi_{\g\oplus\g})=\phi_\g^L-\phi_\g^R=0$, $(G,\rho,0)$ is a quasi-Poisson $\g\oplus\g$-manifold, and it is easily seen to be quasi-symplectic. 
If $G$ is 1-connected, it follows that the pair groupoid
$$(G\times G,\rho_1,\pi):=(G,\rho,0)\circledast(G,\rho,0)$$
is the source-1-connected Hamiltonian quasi-symplectic $\g$-groupoid integrating $(G,\rho,0)$; its moment map is given by $\Phi:(a,b)\to(a\cdot b^{-1},a^{-1}\cdot b)$. 
In \cite{Alekseev00} this example is called the double and is denoted $D(G)$.

\subsection{The group $G$}
The simplest example of a quasi-Poisson $\g$-bialgebra is $\g$ with $\rho=\id$ and $\DIF=0$. We get that $\mu_\rho:\g^*\to\g$ is the identification of $\g^*$ with $\g$ via the inner product and it is an isomorphism of Lie algebras (more on quasi-Poisson $\g$-bialgebras is in \S~\ref{sec:fex}).

The corresponding Hamiltonian quasi-Poisson $\g$-group is $G$ with
\[ \pi_G=\frac{1}{2}\sum_i{e^i}^L\wedge{e_i}^R , \]
$\Phi=\id$, and $\rho$ is the conjugation. It appears in \cite{Alekseev00} as the basic example of a Hamiltonian quasi-Poisson space.


\subsection{Fused double}

Since $D(G)$ is a Hamiltonian quasi-Poisson $\g\oplus\g$-manifold, it has a fusion \eqref{fusion}, which we denote by $\mathbf{D}(G)$ \cite{Alekseev00}. There is a Hamiltonian quasi-symplectic groupoid structure on $\mathbf{D}(G)$ given as follows. The source and target maps, $s,t:\mathbf{D}(G)\rightrightarrows G$, are $s(a,b)=a b^{-1}$, $t(a,b)=b^{-1} a$, and the composition is
\[(a',b')\cdot(a,b)=(a'a,a'b)=(a'a,b'a)\, .\]

This Hamiltonian quasi-symplectic $G$-groupoid is called the AMM groupoid in \cite{Behrend03,Xu03}. It integrates the quasi-Poisson group $G$.

\begin{remark}[Poisson Parallel]
The double $D(G)$ is in some sense the quasi-Poisson analogue of the symplectic groupoid $T^*G$, where $G$ is a Lie group. $T^*G$ integrates the trivial Poisson manifold $G$. It has another groupoid structure, as the groupoid integrating $\g^*$; this is analogous to the groupoid $\mathbf{D}(G)$.
\end{remark}

\begin{remark}
Geometrically, $D(G)$ is the moduli space of flat $\g$-connections on a cylinder, with a marked point on each boundary circle. The composition in $D(G)$ corresponds to cutting cylinders along curves connecting the marked points (the first cylinder along a straight line and the second along a curve that goes once around) and gluing them to form a single cylinder. The composition in $\mathbf{D}(G)$ is just concatenation of cylinders. These two compositions don't form a double groupoid, rather they commute modulo a Dehn twist of the cylinder.
\end{remark}

\subsection{Actions with coisotropic stabilizers}

Let $\mf{q}\subset \g$ be a subalgebra which is coisotropic with respect to the quadratic form, let $\h=\mf{q}^\perp$. Suppose that $\h$ and $\mf{q}$ integrate to closed subgroups $H,Q\subset G$. Let $G$ act on the left of $G/Q$, $\rho':\g\to \mathbf{\Gamma}(T(G/Q))$ be the corresponding map of Lie algebras, and $(G/Q,\rho',0)$ be the quasi-Poisson $\g$-manifold of Example~\ref{coisotropicqpm}. We will be interested in calculating the Hamiltonian quasi-Poisson $\g$-groupoid that $(G/Q,\rho',0)$ integrates to.
Notice that the Lie algebroid $T^*(G/Q)$ has vanishing anchor, i.e.\ it is a bundle of Lie algebras. Consequently, the groupoid is a bundle of groups.

Since $G\times G$ acts on $D(G)$ by groupoid morphisms, we have the following morphism of groupoids 
$$\xymatrix{
D(G)\ar@{}|=[r]&G\times G\ar[r] \ar@<1ex>^s[d]\ar@<-1ex>_t[d]& D(G)/(\{1\}\times Q)\ar@<1ex>[d]\ar@<-1ex>[d]\\
&G\ar[r] & G/Q\\
}$$

By \cite[\S~6]{Alekseev00}, $D(G)/(\{1\}\times Q)$ is a Hamiltonian quasi-Poisson $\g$-manifold, with moment map the first component of $\Phi$; and it follows that \begin{equation}\label{hqpg1}D(G)/(\{1\}\times Q)\rightrightarrows G/Q\end{equation} is a Hamiltonian quasi-Poisson $\g$-groupoid. 

It is easy to check that $$X=\Phi^{-1}(G\times H)/(\{1\}\times Q)=\{(a,b)\in (G\times G)/\on{diag}(Q)\mid a^{-1}\cdot b\in H\}$$  is a subgroupoid of $D(G)/(\{1\}\times Q)$ (where $\on{diag}(Q)\subset G\times G$ denotes the diagonal embedding). However, one may also check that $X$ is precisely the leaf of the foliation corresponding to the Lie algebroid $T^*(D(G)/(\{1\}\times Q))$ which passes through the set of identity elements. Consequently, since \eqref{hqpg1} is a Hamiltonian quasi-Poisson $\g$-groupoid,  $$X\rightrightarrows G/Q$$ is a quasi-Hamiltonian $\g$-groupoid integrating $(G/Q,\rho',0)$.

\begin{remark}[Poisson Parallel]
This construction is a quasi-Poisson analogue of the (rough) principle in Poisson geometry that ``symplectization commutes with reduction'', \cite{Fernandes07,Stefanini07,Fernandes09}.
\end{remark}

\begin{remark}
One can notice that the groupoid $X\rightrightarrows G/Q$ is braided-commutative with respect to the braiding on the category of Hamiltonian quasi-symplectic spaces. More generally, braided-commutative Hamiltonian quasi-symplectic groupoids integrate the quasi-Poisson manifolds of the form $(M,\rho,0)$ where $\rho$ has coisotropic stabilizers. This corresponds to the following fact in Drinfeld's category \cite{Drinfeld90} of modules of $\g$ with braided monoidal structure given by a choice of an associator: braided-commutative algebras in this category are $\g$-modules $A$ with a commutative $\g$-equivariant algebra structure $A\otimes A\to A$, such that $\sum_i(e^i\cdot x)(e_i\cdot y)=0$ for every $x,y\in A$. If $A=C^\infty(M)$, it means exactly that the action of $\g$ on $M$ has coisotropic stabilizers.

Thus one can say in this case that $C^\infty(M)$, with its original product and considered as an object of Drinfeld's category, is the quantization of $(M,\rho,0)$.
\end{remark}


\section{Courant algebroids and Manin pairs}\label{sect:CAandMp}
\subsection{Definition of a Manin pair}
Dirac structures were introduced by Courant-Weinstein in \cite{BeyondPoissonS} (see also \cite{Courant90}) in order to provide a unified setting in which to study closed 2-forms, Poisson structures, and their corresponding Hamiltonian vector fields. Courant algebroids were introduced by Liu-Weinstein-Xu \cite{ManinTriplesBi} to provide an abstract setting from which to study Dirac structures.

\begin{definition}
A {\em Courant algebroid} is a quadruple, $(\mathbb{E},\langle\cdot,\cdot\rangle,\mathbf{a},\lb\cdot,\cdot\rb)$, consisting of a pseudo-euclidean vector bundle $(\mathbb{E}\to M,\langle\cdot,\cdot\rangle)$, a bundle map $\mathbf{a}:\mathbb{E}\to TM$ called the \emph{anchor}, and a bilinear bracket $\lb\cdot,\cdot\rb$ on the space of sections $\mathbf{\Gamma}(\mathbb{E})$ called the \emph{Courant bracket}, such that the following axioms hold
\begin{enumerate}
\renewcommand{\labelenumi}{C-\arabic{enumi}}
\item $\lb X_1,\lb X_2, X_3\rb\rb=\lb\lb X_1,X_2\rb,X_3\rb+\lb X_2,\lb X_1,X_3\rb\rb$\label{C1}
\item $\mathbf{a}(X_1)\langle X_2,X_3\rangle = \langle \lb X_1,X_2\rb,X_3\rangle+\langle X_2,\lb X_1,X_3\rb\rangle$\label{C2}
\item $\lb X_1,X_2\rb+\lb X_2,X_1\rb = \mathbf{a}^*(d\langle X_1,X_2\rangle),$\label{C3}
\end{enumerate}
for $X_i\in\mathbf{\Gamma}(\mathbb{E})$. Here $\mathbf{a}^*:T^*M\to \mathbb{E}^*\cong \mathbb{E}$ is dual to $\mathbf{a}$ using the isomorphism given by inner product $\langle\cdot,\cdot\rangle$. We will often refer to $\mathbb{E}$ as a Courant algebroid, the quadruple $(\mathbb{E},\langle\cdot,\cdot\rangle,\mathbf{a},\lb\cdot,\cdot\rb)$ being understood.

Note that $(\mathbb{E},-\langle\cdot,\cdot\rangle,\mathbf{a},\lb\cdot,\cdot\rb)$ is also a Courant algebroid, which we denote by $\bar{\mathbb{E}}$.

A subbundle $A\subset \mathbb{E}$ of a Courant algebroid is called \emph{Lagrangian} if $A^\perp=A$, and is called a \emph{Dirac structure} if the space of sections $\mathbf{\Gamma}{A}$ is closed under the Courant bracket. A \emph{Manin Pair} is a pair $(\mathbb{E},A)$, where $\mathbb{E}$ is a Courant algebroid, and $A\subset \mathbb{E}$ is a Dirac structure.
\end{definition}

\begin{example}[The standard Courant algebroid]
Let $M$ be a manifold, $\eta\in\Omega^3(M)$ be a closed 3-form, and let $\TT_\eta M=T^*M\oplus TM$ be the pseudo-euclidean vector bundle with the inner product given by the canonical pairing,
$$\langle \alpha+X,\beta+Y\rangle=\alpha(Y)+\beta(X),$$
and the Courant bracket given by
\begin{equation}\label{standbrac}\lb \alpha+X,\beta+Y\rb=[X,Y]+\Lied_X\beta-\iota_Yd\alpha +\iota_Y\iota_X\eta,\end{equation}
for any $\alpha,\beta\in\Omega^*(M)$ and $X,Y\in\mathbf{\Gamma}(TM)$. With the anchor $\mathbf{a}:\TT_\eta M=T^*M\oplus TM\to TM$ defined as the projection along $T^*M$, $\TT_\eta M$ becomes a Courant algebroid, called the \emph{standard Courant algebroid} twisted by $\eta$. We will often write $\TT M:= \TT_0 M$.

Examples of Manin pairs are $(\TT M,TM)$ and $(\TT M,T^*M)$.
\end{example}

\begin{example}[quadratic Lie algebras]\label{ddiagg}
A Courant algebroid over a point is just a quadratic Lie algebra.

Let $\g$ be a quadratic Lie algebra, $\dd=\g\oplus\bar\g$ and $\on{diag}:\g\to\dd$ be the diagonal embedding. Then $(\dd,\on{diag}(\g))$ is a Manin pair.
\end{example}

Suppose $\rho:M\times \dd\to TM$ defines an action of a quadratic Lie algebra on a manifold $M$ with coisotropic stabilizers (see Example~\ref{coisotropicqpm}). Then as shown in \cite{Li-Bland08}, there is a unique Courant algebroid stucture on $M\times\dd$ such that
\begin{itemize}
\item the anchor is $\rho$,
\item the Courant bracket extends the Lie bracket on constant sections, and
\item the pseudo-euclidean structure is given by the quadratic form on $\dd$.
\end{itemize}

\begin{example}[The Cartan-Dirac structure]\label{Ag}\label{CartanDirac}
Let $G$ be a Lie group integrating a quadratic Lie algebra $\g$. Then the action of $\dd=\g\oplus\bar\g$ on $G$ by $$\rho:(g,\xi,\eta)\to-\xi^R_g+\eta^L_g$$ has coisotropic stabilizers. Therefore $\A_G:=G\times\dd$ is a Courant algebroid and $E_G:=G\times\on{diag}(\g)$ is a Dirac structure called the Cartan-Dirac structure \cite{Bursztyn03-1,Kotov-Schaller-Stobl05,Severa-Weinstein01}. It was introduced independently by Alekseev, Strobl and \u{S}evera (see also \cite{PureSpinorsOnL,Li-Bland08}).



\end{example}

\begin{example}\label{ex:s-val}
Suppose that $D$ is a Lie group with quadratic Lie algebra $\dd$ and $G\subset D$ is a closed subgroup with Lagrangian Lie subalgebra $\g\subset \dd$.
Then the left action of $\dd$ on $S:=D/G$ has coisotropic stabilizers. Therefore $(S\times\dd,S\times\g)$ is a Manin pair.
This fact goes back to the unpublished work of the second author \cite{LetToWein}, and A. Alekseev and P. Xu \cite{AlekseevXuDe}. See also \cite{Bursztyn07-1,Ponte08,Bursztyn08,Li-Bland08}.
\end{example}

\subsection{Morphisms of Manin pairs}
To describe morphisms of Manin Pairs, we first need to recall the notions of generalized Dirac structures and Courant morphisms, both due to the second author \cite{LetToWein,Severa01} (see also \cite{AlekseevXuDe,Bursztyn08,Ponte08}).

\begin{definition}[Generalized Dirac structure with support]
Let $\mathbb{E}\to M$ be a Courant algebroid, and $S\subset M$ be a submanifold. Let $\mathbb{E}\rvert_S$ denote the restriction of the pseudo-euclidean vector bundle $\mathbb{E}$ to $S$. A \emph{generalized Dirac structure with support} on $S$ is a subbundle $K\subset \mathbb{E}\rvert_S$ such that
\begin{enumerate}
\renewcommand{\labelenumi}{GD-\arabic{enumi}}
\item $K$ is Lagrangian, namely $K^\perp=K$,
\item $\mathbf{a}(K)\subset TS$, and
\item if $X_i\in\mathbf{\Gamma}(\mathbb{E})$, and $X_i\rvert_S\in\mathbf{\Gamma}(K)$, then $\lb X_1,X_2\rb\rvert_S\in \mathbf{\Gamma}(K)$.
\end{enumerate}
\end{definition}

For any smooth map $f:M\to N$, we let $\on{Gr}_f=\{(m,f(m))\mid m\in M\}$ denote its graph.

\begin{definition}[Courant morphism] If $\mathbb{E}\to M$ and $\mathbb{F}\to M$ denote two Courant algebroids, a \emph{Courant morphism} $(f,K):\mathbb{E}\dasharrow\mathbb{F}$ is a smooth map $f:M\to N$ together with a generalized Dirac structure $K\subset\mathbb{F}\times\bar{\mathbb{E}}$ with support on $\on{Gr}_f$, the graph of $f$.
\end{definition}

\begin{definition}[Morphism of Manin Pairs]\label{morphmanpair}
Suppose $(\mathbb{E},A)$ and $(\mathbb{F},B)$ are Manin pairs. A \emph{morphism of Manin pairs} $(f,K):(\mathbb{E},A)\dasharrow(\mathbb{F},B)$ is a Courant morphism $(f,K):\mathbb{E}\dasharrow\mathbb{F}$ such that the image of $K$ under the projection $\mathbb{F}\times\bar{\mathbb{E}}\to\mathbb{F}/B\times\bar{\mathbb{E}}/A$ is the graph of a bundle map $$\phi_K:A^*\to f^*B^*.$$ Here we used $A^*\cong \mathbb{E}/A$, and $B^*\cong \mathbb{F}/B$.
\end{definition}

The notion of a \emph{Morphism of Manin Pairs} was introduced in \cite{Bursztyn08} to study general moment maps. The following example is also found in \cite{Bursztyn08}.



\begin{example}[Strong Dirac Morphisms]\label{strongexamp}
Let $f:M\to N$ be a map between smooth manifolds, let $\xi\in\Omega^3(M)$ and $\eta\in\Omega^3(N)$ be closed 3-forms, and $\omega\in\Omega^2(M)$ a 2-form. Then
$$K_{(f,\omega)}=\{(f^*\alpha-\iota_X\omega,X,\alpha,f_*X)\mid p\in M,\alpha\in T^*_{f(p)}N, X\in T_pM\}$$
is a generalized Dirac structure of $\TT_\xi M\times\bar\TT_\eta N$ supported on $\on{Gr}_f$ if and only if $\xi=f^*\eta+d\omega$. 

Suppose now $(\TT_\xi M,A)$ and $(\TT_\eta N,B)$ are two Manin pairs. $(f,K_{(f,\omega)})$ define a morphism of Manin pairs 
$$(f,K_{(f,\omega)}):(\TT_\xi M,A)\dasharrow(\TT_\eta N,B)$$
if and only if $(f,\omega):(M,A,\xi)\to (N,B,\eta)$ is a strong Dirac morphism, as in \cite{PureSpinorsOnL}. With $\omega=0$, $(f,K_{(f,0)}):(\TT_{f^*\eta}M,A)\dasharrow(\TT_\eta N,B)$ is a morphism of Manin pairs if and only if $f$ is a strong Dirac map from $A$ to $B$, as in \cite{Bursztyn03,Bursztyn07-1}.
\end{example}

In order to identify morphisms of Manin Pairs of the form given in Example~\ref{strongexamp} we introduce the notion of

\begin{definition}[Full Morphisms of Manin Pairs]\label{fullmorph}
A morphism of Manin pairs, $(f,K):(\mathbb{E},A)\to(\mathbb{F},B)$, is called \emph{full} if $\mathbf{a}(K)=T\on{Gr}_f$.
\end{definition}

\begin{remark}\label{strongdiracremark}
Suppose $(f,K):(\TT_\xi M,A)\dasharrow(\TT_\eta N,B)$ is a \emph{full} morphism of Manin pairs. Since $\mathbf{a}:K\to T\on{Gr}_f$ is surjective, there must be a 2-form $\omega\in\Omega^2(M)$ such that $K=K_{(f,\omega)}$. Furthermore, the morphism of Manin pairs must be of the form given in Example~\ref{strongexamp}. So $(f,K)$ simply describes a strong Dirac morphism $(f,\omega):(M,A,\xi)\to(N,B,\eta)$.
\end{remark}

\subsection{Multiplicative Manin pairs}
We will be interested in Dirac structures living on groupoids which are multiplicative in some sense. Although Poisson Lie groups \cite{Drinfeld83,Semenov-Tian-Shansky85,lu90,thesis-3}, Poisson groupoids \cite{weinstein87,Mackenzie-Xu94} and symplectic groupoids \cite{weinstein87-1,SymplGrpoids(W} are examples of such objects, the first comprehensive study of them appears in the papers of Ortiz \cite{Ortiz08,(MultDiracSt} (see also \cite{Bursztyn09}), where \emph{multiplicative Dirac structures} are defined. In this section we show that, when endowed with certain morphisms, multiplicative Dirac structures form subcategory of the category of Manin pairs. This subcategory, called \emph{Multiplicative Manin Pairs}, should be thought of as the groupoid objects in the category of Manin Pairs.

A Lie groupoid $V\rightrightarrows E$ is called a {\em $\mathcal{VB}$-groupoid}, if $V$ and $E$ are also vector bundles over $\Gamma$ and $M$ respectively, and all the structure maps are smooth vector bundle maps (see Appendix~A, Page~\pageref{vbgroupoids}). In this case $\Gamma\rightrightarrows M$ inherits the structure of a Lie groupoid. We may refer to $V\to \Gamma$ as a $\mathcal{VB}$-groupoid when we want to specify $\Gamma$ as the base of the vector bundle. 

Next we need the concepts of Courant groupoids and multiplicative Dirac structure as given by Mehta \cite{mehta09} and Ortiz \cite{Ortiz08,(MultDiracSt}.

\begin{definition}[Courant groupoid]
Let $\mathbb{E}\to\Gamma$ be a $\mathcal{VB}$-groupoid such that $\mathbb{E}$ is a Courant algebroid, and let
$$K_m\subset\bar{\mathbb{E}}\times\bar{\mathbb{E}}\times\mathbb{E}$$ and 
$$\on{Gr}_{m_\Gamma}\subset\Gamma\times\Gamma\times\Gamma$$
denote the respective graphs of the multiplications for $\mathbb{E}$ and $\Gamma$. $\mathbb{E}$ is called a \emph{Courant groupoid} if $K_m$ is a generalized Dirac structure with support on $\on{Gr}_{m_\Gamma}$.

A \emph{morphism of Courant groupoids} $(f,K):\mathbb{E}\dasharrow\mathbb{F}$ is a Courant morphism for which $f$ is a morphism of Lie groupoids and $K\subset\mathbb{F}\times\bar{\mathbb{E}}$ is a Lie subgroupoid.

A \emph{multiplicative Dirac structure}, is a Dirac structure $A\subset\mathbb{E}$ which is also a subgroupoid. 
A \emph{multiplicative Manin pair} is a pair $(\mathbb{E},A)$, where $\mathbb{E}$ is a Courant groupoid, and $A\subset \mathbb{E}$ is a multiplicative Dirac structure.
\end{definition}

%


\begin{definition}[Morphism of Multiplicative Manin pairs]

A \emph{morphism of multiplicative Manin pairs} $$(f,K):(\mathbb{E},A)\dasharrow(\mathbb{F},B)$$ is a morphism of Manin pairs such that $(f,K):\mathbb{E}\dasharrow\mathbb{F}$ defines a morphism of Courant groupoids.
\end{definition}

The following example is found in \cite{Bursztyn03-1,Ortiz08,Bursztyn09,(MultDiracSt}.
\begin{example}
If $\Gamma$ is any groupoid, then by applying the tangent functor to all the spaces and morphisms $T\Gamma$ becomes a $\mathcal{VB}$-groupoid. $T^*\Gamma$ is its dual $\mathcal{VB}$-groupoid. Consequently $\TT\Gamma=T^*\Gamma\oplus T\Gamma$ becomes a $\mathcal{VB}$-groupoid. One can check that it is actually a Courant groupoid.



Examples of multiplicative Manin pairs are $(\TT\Gamma,T\Gamma)$ and $(\TT\Gamma,T^*\Gamma)$.
\end{example}

\begin{example}\label{ddiaggm}
Let $\g$ be a quadratic Lie algebra, and consider $\dd=\g\oplus\bar\g$ together with the pair groupoid structure. $(\dd,\on{diag}(\g))$ is a multiplicative Manin pair.
\end{example}

The following example is found in \cite{Bursztyn03,PureSpinorsOnL,Li-Bland08}.
\begin{example}\label{AgEg}
The Courant algebroid $\A_G=G\times\g\oplus\bar\g$ (Definition~\ref{Ag}) is the direct product of the Lie group $G$ and the pair groupoid $\g\oplus\bar\g$; with this groupoid structure, it becomes a Courant groupoid. Clearly the Cartan Dirac structure $E_G$ (Definition~\ref{CartanDirac}) is a subgroupoid, and so $(\A_G,E_G)$ is a multiplicative Manin pair.
\end{example}






The infinitesimal version of a multiplicative Dirac structure is studied in \cite{(MultDiracSt} (see also \cite{Bursztyn03-1,Ortiz08,Bursztyn09}). We will be interested in the infinitesimal version of {\em morphisms} of multiplicative Manin pairs. To study this in \S~\ref{sect:MpandMPmflds}, we will find it more convenient to use an alternative description of Manin pairs via graded Poisson geometry \cite{Bursztyn08}.
%


\section{Manin pairs and MP-manifolds}\label{sect:MpandMPmflds}
It was shown in \cite{Bursztyn08} that category of Manin pairs is equivalent to the category of MP-manifolds (first introduced in \cite{LetToWein}, see also \cite{NonComDiffForm,PoissonActions}). Since using the latter category can sometimes bring more geometric insight, we recall the equivalence described in \cite{Bursztyn08}.

\begin{definition}
An {\em MP-manifold}  is a principal $\mathbb{R}[2]$ bundle $P\to A^*[1]$, where $A$ is a vector bundle over a manifold $M$, such that $P$ carries a Poisson structure of degree $-1$ which is $\mathbb{R}[2]$ invariant. We call $M$ the MP-base of $P$.
\end{definition}

There is a natural notion of morphisms for MP-manifolds:

\begin{definition}[Morphisms of MP-manifolds]\label{MPmorphdef}
Let $P\to A^*[1]$ and $Q\to B^*[1]$ be two MP-manifolds, then a morphism of MP-manifolds

$$\xymatrix{
P \ar[d] \ar[r]^F & Q\ar[d]\\
A^*[1] \ar[r]^{F/\mathbb{R}[2]}&B^*[1]\\
}$$
is an $\mathbb{R}[2]$ equivariant Poisson map. We will often abbreviate morphisms of MP-manifolds as $F:P\dasharrow Q$.
\end{definition}

The Poisson structure on $P$ induces a map $\pi^\sharp:T^*[2]P\to T[1]P$. If we use $F$ to pull back the cotangent bundle on $Q$ to a bundle over $P$, then we have a map $(T[-2]F)^*:F^*(T^*[2]Q)\to T^*[2]P$ dual to the tangent map. Finally, if $M$ is the MP-base of $P$ we denote the projection $p:P\to A^*[1]\to M$, and the corresponding tangent map $T[1]p:T[1]P\to T[1]M$.

\begin{definition}
A morphism of MP-manifolds is called \emph{full} if the map
\begin{equation}\label{fullcond}T[1]p\circ\pi^\sharp\circ (T[-2]F)^*:F^*(T^*[2]Q)\to T[1]M\end{equation}
is surjective.
\end{definition}

\begin{remark}\label{immtoinj}
Let $\mu:T^*[2]Q\to \mathbb{R}$ be the moment map for the $\mathbb{R}[2]$ action. Let $F^*(\mu^{-1}(1))$ denote the pullback of $\mu^{-1}(1)$ by $F$. Since $\pi$ is $\mathbb{R}[2]$ invariant and $F$ is $\mathbb{R}[2]$ equivariant, \eqref{fullcond} is surjective if and only if
\begin{equation}\label{fullcond2}(T[1]p/\mathbb{R}[2])\circ (\pi^\sharp/\mathbb{R}[2])\circ ((T[-2]F)^*/\mathbb{R}[2]):F^*(\mu^{-1}(1))/\mathbb{R}[2]\to T[1]M\end{equation}
is surjective.

However the $1$-graded part of $F^*(\mu^{-1}(1))/\mathbb{R}[2]$ describes a vector bundle $K_F\to M$, so \eqref{fullcond2} describes a base-preserving morphism of vector bundles $\mathbf{a}_F:K_F\to TM$. It is clear that \eqref{fullcond2} is surjective if and only if $\mathbf{a}_F$ is surjective, or equivalently 
\begin{equation}\label{fullcond3}\mathbf{a}_F^*:T^*M\to K^*_F\end{equation}
is injective. Note that since $\mathbf{a}^*_F$ is a smooth base preserving morphism of vector bundles, it is injective if and only if it is an immersion.
\end{remark}

\begin{theorem}
The equivalence between the categories of Manin pairs and MP-manifolds described in \cite{Bursztyn08} identifies full morphisms of Manin pairs with full morphisms of MP-manifolds.
\end{theorem}
\begin{proof}
First we recall the functor from the category of MP-manifolds to Manin pairs given in \cite{Bursztyn08}.

In the work of Roytenberg, Vaintrob, Weinstein and the second author \cite{LetToWein,Roytenberg99,Severa01,Roytenberg02, GradedSymplSup}, it was shown that a Courant algebroid $\mathbb{E}$ is equivalent to a degree 3 function, $\Theta$, on a non-negatively graded degree 2 symplectic manifold, $\mathcal{M}$, such that $\{\Theta,\Theta\}=0$. In this picture, a Dirac structure with support on a submanifold corresponds to a Lagrangian submanifold of $\mathcal{M}$ on which $\Theta$ vanishes.

Suppose $P$ is an MP-manifold. The Poisson structure on $P$ corresponds to an $\mathbb{R}[2]$-invariant degree three function, $\pi\in C^\infty(T^*[2]P)$, such that $\{\pi,\pi\}=0$. $\pi$ descends to a degree 3 function, $\Theta$, on $\mathcal{M}_P=T^*[2]P/\!\!/_1\mathbb{R}[2]$, the symplectic reduction at moment value 1, such that $\{\Theta,\Theta\}=0$. Thus $P$ defines a Courant algebroid, $\mathbb{E}$. The map $\mathcal{M}_P\to A^*[1]$ corresponds to a base-preserving vector bundle morphism
\begin{equation}\label{projection}\mathbb{E}\to A^*\end{equation}
whose kernel is a Dirac structure $A\subset \mathbb{E}$. In this way, $P$ defines a Manin pair $(\mathbb{E},A)$.

Suppose that $P$ and $Q$ are MP-manifolds corresponding to the Manin pairs $(\mathbb{E},A)$ and $(\mathbb{F},B)$. Let $F:P\dasharrow Q$ be a morphism of MP-manifolds, and $\on{Gr}_F$ its graph. Let $M$ and $N$ be the respective MP-bases, and $f:M\to N$ denote the restriction of $F$. The conormal bundle $N^*[2]\on{Gr}_F\subset T^*[2]Q\times\overline{T^*[2]P}$ is a Lagrangian submanifold on which $\pi_{Q\times\bar P}$ vanishes. The reduction $L_F$ of $N^*[2]\on{Gr}_F$ to the symplectic quotient $\mathcal{M}_Q\times\overline{\mathcal{M_P}}$ is a Lagrangian submanifold on which $\Theta_{Q\times\bar P}$ vanishes; it corresponds to a generalized Dirac structure $K$ defining a morphism of Manin pairs $(f,K):(\mathbb{E},A)\dasharrow(\mathbb{F},B)$ (see \cite{Bursztyn08} for details).

We need to show that $(f,K)$ is full if and only if $F$ is full. By identifying $M$ with $\on{Gr}_f$, the anchor map takes the form $\mathbf{a}:K\to TM$. It is surjective if and only if $(f,K)$ is full. By identifying $P$ with $\on{Gr}_F$, the conormal bundle $N^*[2]\on{Gr}_F$ may be naturally identified with $F^*(T^*[2]Q)$, and similarly $L_F$ with $F^*(\mu^{-1}(1))/\mathbb{R}[2]$ (where $\mu:T^*[2]Q\to\mathbb{R}$ is the moment map for the $\mathbb{R}[2]$ action, as in Remark~\ref{immtoinj}). Under this identification, the vector bundle $K_F$ of Remark~\ref{immtoinj} corresponds to $K$, and $\mathbf{a}_F$ to the anchor map $\mathbf{a}$. Consequently by Remark~\ref{immtoinj}, $(f,K)$ is a full morphism of Manin pairs if and only if $F$ is a full morphism of MP-manifolds. 

\end{proof}

\begin{example}\label{T^*[1]M}
The simplest example of an MP-manifold is $T^*[1]M\times\mathbb{R}[2]$, where the Poisson structure comes from the canonical symplectic structure on the cotangent bundle and the trivial one on $\mathbb{R}[2]$, and the $\mathbb{R}[2]$ action is the obvious one. It corresponds to the Manin pair $(\TT M, TM)$.
\end{example}

\subsection{Multiplicative Manin pairs and MP-groupoids}
In this section we introduce the category of MP-groupoids (a subcategory of MP-manifolds) and establish an equivalence between it and the category multiplicative Manin pairs.

Given the notion of a Poisson groupoid \cite{weinstein87,Mackenzie-Xu94}, there is a natural notion of MP-groupoids.

\begin{definition}[MP-groupoid]\label{MPG}
An MP-manifold $P\to A^*[1]$ is called an \emph{MP-groupoid}, if it is a Poisson groupoid, and the $\mathbb{R}[2]$ action map, $P\times\mathbb{R}[2]\to P$, is a groupoid morphism.

In more detail, let $\bar P$ be the graded groupoid $P$ with minus the Poisson structure \eqref{barM}. Let
$$\on{Gr}_{m_P}\subset P\times P\times \bar P$$
denote the graph of the multiplication and $P_0\subset P$ the submanifold of identity elements. Then $P$ is a MP-groupoid if $\on{Gr}_{m_P}$ and $P_0$ are coisotropic submanifolds, and the action map $P\times\mathbb{R}[2]\to P$ is a groupoid morphism.

A morphism of MP-groupoids is a morphism of MP-manifolds which is also a groupoid morphism.

If an MP-groupoid is actually a (graded) Lie group, we may refer to it as an MP-group.
\end{definition}

The following proposition should come as no surprise.

\begin{proposition}\label{MPgmMPe}
The equivalence between the category of Manin pairs and MP-manifolds induces an equivalence between the categories of multiplicative Manin Pairs and MP-groupoids.
\end{proposition}

First we need a short lemma.

\begin{lemma}\label{multequiv}
Let $\mathbb{E}\to\Gamma$ be a Courant groupoid. A Dirac structure $A\subset\mathbb{E}$ is multiplicative if and only if there is a $\mathcal{VB}$-groupoid structure on $A^*\to\Gamma$ for which the projection $p:\mathbb{E}\to\mathbb{E}/A\cong A^*$ is a morphism of $\mathcal{VB}$-groupoids.

Furthermore, if $(f,K):(\mathbb{E},A)\dasharrow(\mathbb{F},B)$ is a morphism of multiplicative Manin pairs, then the induced map $\phi_K:A^*\to B^*$ is a morphism of Lie groupoids. 
\end{lemma}
\begin{proof}
We recall \cite[Proposition 11.2.5]{Mackenzie05} that $A^*$ has the structure of a $\mathcal{VB}$-groupoid if and only if $A$ also has the structure of a $\mathcal{VB}$-groupoid (See Appendix~A, Page~\pageref{vbgroupoids} for details). Furthermore the inner product $\langle\cdot,\cdot\rangle$ defines an isomorphism between the $\mathcal{VB}$-groupoids $\mathbb{E}$ and $\mathbb{E}^*$, where the latter has the structure of the dual $\mathcal{VB}$-groupoid. Using this isomorphism, the projection $p:\mathbb{E}\to A^*$ is dual to the inclusion $i:A\to\mathbb{E}$. \cite[Proposition 11.2.6]{Mackenzie05} states that if either $p$ or $i$ is a morphism of $\mathcal{VB}$-groupoids, then they both are.

Since $K\subset \mathbb{E}\times\overline{\mathbb{F}}$ is a subgroupoid, it follows that $K/\big(K\cap (A\times B)\big)$ is a subgroupoid of $(\mathbb{E}\times\overline{\mathbb{F}})/(A\times B)=A^*\times B^*$. However, $K/\big(K\cap (A\times B)\big)$ is the graph of $\phi_K:A^*\to B^*$. Therefore $\phi_K$ is a morphism of Lie groupoids.
\end{proof}

\begin{proof}[Sketch of proof of Proposition~\ref{MPgmMPe}]
Let $P$ be an MP-groupoid. $\mathcal{M}=T^*[2]P/\!\!/_1\mathbb{R}[2]$ is a symplectic groupoid, and the degree 3 function on $\mathcal{M}$ corresponding to the Poisson structure on $P$ is multiplicative. In other words $\mathcal{M}$ describes a Courant groupoid $\mathbb{E}$. Furthermore, the map \eqref{projection} describes a morphism of $\mathcal{VB}$-groupoids. Therefore, by Lemma~\ref{multequiv}, the Manin pair $(\mathbb{E},A)$ corresponding to $P$ is multiplicative. Similarly, it is easy to check that morphisms of MP-groupoids correspond to morphisms of multiplicative Manin pairs.
\end{proof}

\subsection{MP-algebroids}

We will be interested in studying the infinitesimal counterparts of MP-groupoids. Intuitively, since an MP-groupoid is just a Poisson groupoid (together with some free $\mathbb{R}[2]$ action), it must integrate some Lie bialgebroid \cite{Mackenzie-Xu94} (and since the $\mathbb{R}[2]$ action is given by a groupoid morphism, it must integrate a Lie algebroid morphism). This intuition should motivate the following definition.

\begin{definition}[MP-algebroid]\label{MPA}
An MP-algebroid is a graded Lie algebroid $P$, such that $P$ is also an MP-manifold, and
\begin{enumerate}
\renewcommand{\labelenumi}{MPA-\arabic{enumi}}
\item the Poisson structure on $P$ is linear, defining a Lie algebroid structure on $P^*$ (see \cite{LieAlgebroidsH}), 
\item the Lie algebroid structures on $P$ and $P^*$ are compatible, so that $P$ is a Lie bialgebroid (see \cite{Mackenzie-Xu94,Voronov02,Voronov06,Voronov07}), and
\item the action map $P\times\mathbb{R}[2]\to P$ is a Lie algebroid morphism, where $\mathbb{R}[2]$ is viewed as a trivial Lie algebra.
\end{enumerate}

We call $P$ integrable if it is integrable as a Lie algebroid, and if $P$ is actually a Lie algebra (rather than just a Lie algebroid), we may call it an MP-algebra.

Morphisms of MP-algebroids are morphisms of Lie algebroids which are also morphisms of MP-manifolds.
\end{definition}

\begin{proposition}
The category of integrable MP-algebroids is equivalent to the category of source 1-connected MP-groupoids.
\end{proposition}
\begin{proof}
The proofs in Mackenzie-Xu \cite{Mackenzie-Xu94} apply in the graded category to establish an equivalence between the categories of source 1-connected graded Poisson groupoids and integrable graded Lie bialgebroids. Since MP-groupoids and MP-algebroids are simply Poisson groupoids and Lie bialgebroids (respectively) together with additional requirements regarding an $\mathbb{R}[2]$ action, we need only show that these requirements correspond to each other under the Mackenzie-Xu equivalence.

We may view $\mathbb{R}[2]$ as a trivial Lie bialgebroid, which integrates to $\mathbb{R}[2]$ (viewed as a groupoid under addition with the trivial Poisson structure). On the Poisson groupoid level we required the existence of an $\mathbb{R}[2]$ action map which was both a groupoid and a Poisson morphism (see Definition~\ref{MPG}), clearly this corresponds on the Lie bialgebroid level to requiring the existence of an $\mathbb{R}[2]$ action map which is both a Lie algebroid and a Poisson morphism (see Definition~\ref{MPA}).
\end{proof}

\begin{example}\label{T*A}
Let $B$ be a Lie algebroid integrating to the groupoid $\Gamma$. Then $T^*[1]B\times\mathbb{R}[2]$  is a MP-algebroid which integrates to the MP-groupoid $T^*[1]\Gamma\times\mathbb{R}[2]$, where the cotangent bundle has the canonical Poisson structure and $\mathbb{R}[2]$ has the trivial one. This MP-groupoid corresponds to the multiplicative Manin pair $(\TT\Gamma,T\Gamma)$.
\end{example}

\begin{theorem}\label{inffull}
Let $P$ and $Q$ be MP-algebroids integrating to MP-groupoids $\Gamma_P$ and $\Gamma_Q$. A morphism of MP-groupoids $F:\Gamma_P\dasharrow\Gamma_Q$ is full if and only if the corresponding morphism of MP-algebroids $f:P\dasharrow Q$ is full.
\end{theorem}

\begin{proof}

Let $\Gamma$ be the MP-base of $\Gamma_P$ and let the Lie algebroid of $\Gamma$ be denoted by $A$. Clearly $A$ is the MP-base of $P$.

$F$ is full if and only if the morphism 
\begin{subequations}
\begin{equation}\label{fullgrp}\mathbf{a}^*_F:T^*\Gamma\to K_F\end{equation}
described in \eqref{fullcond3} is an injective immersion. Meanwhile $f$ is full if and only if 
\begin{equation}\label{fullalg}\mathbf{a}^*_f:T^*A\to K_f\end{equation}\end{subequations}
is an injective immersion.

So $f$ is full if and only if \eqref{fullalg} describes the inclusion of a subalgebroid, while $F$ is full if and only if \eqref{fullgrp} describes the inclusion of a subgroupoid. However \eqref{fullgrp} integrates the Lie algebroid morphism \eqref{fullalg}. Consequently, if $F$ is full, then so is $f$. On the other hand, if $f$ is full, then \eqref{fullgrp} is an immersion (by \cite[\S~3.2]{Moerdijk02}), and consequently $F$ is full (by Remark~\ref{immtoinj}).
\end{proof}

\subsection{MP groups}
We can now give a description of MP Lie groups in terms of generalized Manin triples. This description is a generalization of the usual description of Lie bialgebras and Poisson-Lie groups. 

Let us remark that an MP group with the base $H$ (where $H$ is a Lie group) is equivalent to a multiplicative Manin pair $(\mathbb{E},A)$ on $H$ such that $\mathbb{E}/A$ is a group (not just a groupoid); equivalently, the space of objects of the groupoid $\mathbb{E}$ is the fiber of $A$ at $1\in H$.

MP groups are the general type of Manin pairs that lead to moment maps admitting a fusion product. If $P$ is an MP group, a $P$-type moment map is a graded Poisson map $T^*[1]M\to P$, and such maps can be multiplied via the product in $P$.

Our generalized Manin triples are $(\mf{f},\mf{h},\mf{k})$, where $\mf{f}$ is a Lie algebra with a chosen ad-invariant element $s$ of $S^2\mf{f}$ and $\mf{h}$, $\mf{k}$ are its subalgebras such that $\mf{f}=\mf{h}\oplus\mf{k}$ as a vector space and $\mf{k}$ is $s$-coisotropic. As we shall see, this data is equivalent to a MP group with the base $H$, the 1-connected group integrating $\mf{h}$.

Let us consider the graded Lie algebra
\[\mc{Q}_s(\mf{f})=\mathbb{R}[2]\oplus\mf{f}^*[1]\oplus\mf{f}\oplus\mathbb{R}[-1]\]
with the Lie bracket given by ($\alpha,\beta\in\mf{f}^*[1]$, $\xi,\eta\in\mf{f}$)
\begin{align*}
[\alpha,\beta]&= s(\alpha,\beta) \cdot T  \\
[\xi,\alpha]&=-\ad(\xi)^*\alpha &[\xi,\eta]&=[\xi,\eta]_{\mf{f}}\\
[D,\alpha]&=s^\sharp(\alpha) &[D,\xi]&=0 &[D,D]&=0
\end{align*}
where $T$ and $D$ are the generators of $\mathbb{R}[2]$ and of $\mathbb{R}[-1]$ respectively, and $T$ is central. It has a non-degenerate pairing of degree 1 given by $\langle T,D\rangle=1$, $\langle \xi,\alpha\rangle=\alpha(\xi)$.

A generalized Manin triple $(\mf{f},\mf{h},\mf{k})$ then gives rise to a pair of transverse Lagrangian subalgebras of $\mc{Q}_s(\mf{f})$
\[\mathbb{R}[2]\oplus\mf{h}^\perp[1]\oplus\mf{h}\quad\text{and}\quad
\mf{k}^\perp[1]\oplus\mf{k}\oplus\mathbb{R}[-1]\ ,\]
i.e.\ to a graded Lie bialgebra (with cobracket of degree $-1$). It makes the graded Poisson-Lie group integrating $\mathbb{R}[2]\oplus\mf{h}^\perp[1]\oplus\mf{h}$ to an MP-group.

\begin{theorem}
An MP group with a 1-connected base $H$ is equivalent to a generalized Manin triple $(\mf{f},\mf{h},\mf{k})$. The corresponding Courant algebroid on $H$ is exact if and only if $s$ is non-degenerate and $\mf{k}\subset\mf{f}$ is Lagrangian.
\end{theorem}

\begin{proof}
MP groups with base $H$ correspond to graded Lie bialgebras (with cobracket $\delta$ of degree $-1$) of the form 
\[\mathbb{R}[2]\oplus V[1]\oplus\mf{h}\]
(where $V$ is some vector space), such that the generator $T$ of $\mathbb{R}[2]$ is central and $\delta(T)=0$. One can easily check that these are exactly the Lie bialgebras coming from triples $(\mf{f},\mf{h},\mf{k})$. 

The Courant algebroid corresponding to the MP-group is transitive if and only if the identity morphism of Manin pairs is full. By Theorem~\ref{inffull} this is equivalent to the identity morphism of MP-algebras being full. That is to say $s^\sharp\rvert_{\mf{k}^\perp}:\mf{k}^\perp\to\mf{k}$ is an injection. 

For dimensional reasons, the Courant algebroid corresponding to the MP-group is exact if and only if $s^\sharp\rvert_{\mf{k}^\perp}:\mf{k}^\perp\to\mf{k}$ is an isomorphism; this means that $s$ is non-degenerate and $\mf{k}\subset\mf{f}$ is Lagrangian.
\end{proof}

\begin{remark}
When the Courant algebroid is exact and moreover the projection of $s$ to $S^2\mf{h}$ is non-degenerate, the corresponding Manin pair on $H$ was used in \cite{KS97} as a boundary condition for the WZW model on the group $H$.
\end{remark}

\begin{remark}
If the Courant algebroid on $H$ is exact then the dual Lie algebra
$\mf{k}^\perp[1]\oplus\mf{k}\oplus\mathbb{R}[-1]$ is isomorphic to $\hat{\mf{k}}$. This type is the most interesting case from the point of view of moment map theory.

  On the other hand, any cobracket $\delta$ of degree $-1$ on $\hat{\mf{k}}$ making $\hat{\mf{k}}$ to a Lie bialgebra comes from a triple 
$(\mf{f},\mf{h},\mf{k})$ with $s$ non-degenerate and $\mf{h}$ Lagrangian.
\end{remark}

\begin{example}
If $s$ is non-degenerate and $\mf{k}\subset\mf{f}$ is Lagrangian, then one may define an invariant element $\eta\in\wedge^3\mf{h}^*$ by $\eta(X,Y,Z)=s^{-1}([X,Y],Z)$ for $X,Y,Z\in\mf{h}$. This corresponds to an invariant closed 3-form on $H$ (the characteristic class of the corresponding exact Courant algebroid). In particular, if $\mf{h}$ is also Lagrangian, then $\eta=0$ and the Courant groupoid over $H$ is just $\TT H$. In this case $\mf{h}$ becomes a Lie bialgebra and the multiplicative Dirac structure $A\subset \TT H$ just describes the corresponding Poisson Lie structure on $H$.

More generally, in \cite{Ortiz08} C. Ortiz classified all multiplicative Dirac structures $A\subset\TT H$ (not simply the ones for which $\TT H/A$ becomes a Lie group).
\end{example}

\section{Manin Pairs and quasi-Poisson structures}\label{sect:MpandqPstruct}
\subsection{Reinterpretation of \S~\ref{part1} in terms of MP-manifolds}\label{atqps}
All the theory described in Part~\ref{part1} can be reinterpreted in terms of MP-manifolds. To begin with,
if $P$ is any MP-manifold, a Poisson map $F:T^*[1]M\to P$ can be canonically lifted to a map of MP-manifolds
$$\tilde F:T^*[1]M\times\mathbb{R}[2]\dasharrow P$$
given by $\tilde F:(x,t)\to F(x)+t$ for any $x\in T^*[1]M$ and $t\in\mathbb{R}[2]$, where the addition refers to the action of $\mathbb{R}[2]$ on $P$. Conversely, given any morphism of MP-manifolds $G:T^*[1]M\times\mathbb{R}[2]\dasharrow P$, the map $F:T^*[1]M\to P$ given by $F:x\to G(x,0)$ is a Poisson morphism; and we may recover $G$ from $F$ since $G=\tilde F$.

We recall from Example~\ref{T^*[1]M}, that $T^*[1]M\times\mathbb{R}[2]$ corresponds to the Manin pair $(\TT M,TM)$. Let $(\mathbb{E},A)$ denote the Manin pair corresponding to $P$. A Poisson morphism $T^*[1]M\to P$ corresponds to a morphism of Manin pairs $(\TT M,TM)\dasharrow (\mathbb{E},A)$.

 The Manin triple $(\mathcal{Q}(\dd),\mathbb{R}[2]\oplus\g[1]\oplus\bar\g,\hat\g)$ of \S~\ref{qhqpgmr} defines a Lie bialgebra structure on $\mathbb{R}[2]\oplus\g[1]\oplus\bar\g$. In other words, $\mathbb{R}[2]\oplus\g[1]\oplus\bar\g$ comes with both a Lie algebra structure and a compatible Poisson bracket of degree $-1$. This together with the natural action of $\mathbb{R}[2]$ gives $\mathbb{R}[2]\oplus\g[1]\oplus\bar\g$ the structure of an MP-algebra. It integrates to the MP-group $\mathbf{G}_{big}$ described in \S~\ref{qhqpgmr}, where $\mathbb{R}[2]$ acts in the obvious way. One may check that the MP-group $\mathbf{G}_{big}$ corresponds to the multiplicative Manin pair $(\A_G,E_G)$ of Example~\ref{AgEg}.

In \S~\ref{qhqpgmr}, we showed that a Hamiltonian quasi-Poisson $\g$-structure on $M$  was equivalent to a Poisson map $T^*[1]M\to \mathbf{G}_{big}$. It is now clear that it also corresponds to a morphism of MP-manifolds $T^*[1]M\times\mathbb{R}[2]\dasharrow \mathbf{G}_{big}$, or simply a morphism of Manin pairs

\begin{equation}\label{qpmpm1}(\TT M,TM)\dasharrow (\A_G,E_G).\end{equation}

This fact was already known to be a direct consequence of \cite[Proposition 3.5]{Bursztyn08} (or of \cite[ Theorem 3.7]{Bursztyn08} and \cite[Theorem 5.22]{PureSpinorsOnL}). As a result of Remark~\ref{strongdiracremark} and \cite[Theorem 5.2]{PureSpinorsOnL}, it is also clear that a Hamiltonian quasi-symplectic $\g$-structure on $M$ corresponds to a full morphism of Manin pairs $(\TT M,TM)\dasharrow (\A_G,E_G)$ or equivalently a full morphism of MP-manifolds $T^*[1]M\times\mathbb{R}[2]\dasharrow \mathbf{G}_{big}$. Furthermore, if $\Gamma$ is a Lie groupoid, then it follows from \S~\ref{qhqpgr} that a compatible Hamiltonian quasi-symplectic $\g$-structure on $\Gamma$ is equivalent to a full morphism of MP-groupoids
$$T^*[1]\Gamma\times\mathbb{R}[2]\dasharrow \mathbf{G}_{big}.$$ The latter is equivalent to a full morphism of multiplicative Manin pairs $$(\TT \Gamma, T\Gamma)\dasharrow (\A_G,E_G).$$

Next consider the Manin triple $(\mathcal{Q}(\g),\mathbb{R}[2]\oplus\g[1],\g\oplus\mathbb{R}[-1])$. The Lie bialgebra structure on $\mathbb{R}[2]\oplus\g[1]$ it defines corresponds to an MP-algebra structure. $\mathbb{R}[2]\oplus\g[1]$ integrates to the MP-group $\mathbf{G}_{small}$ corresponding the the multiplicative Manin pair $(\dd=\g\oplus\bar\g,\on{diag}(\g))$ from Example~\ref{ddiaggm}. The theory in \S~\ref{mathbbR2timesg1} then shows that a quasi-Poisson $\g$ structure on a manifold $M$ corresponds to a morphism of MP-manifolds

$$T^*[1]M\times\mathbb{R}[2]\dasharrow \mathbf{G}_{small}.$$ The latter is equivalent to a morphism of Manin pairs \begin{equation}\label{morphManinPairDesc1}(\TT M,TM)\dasharrow (\dd,\on{diag}(\g)).\end{equation} More generally, Remark~\ref{rem:alt_qPbialg} states that a quasi-Poisson $\g$-bialgebroid structure on a vector bundle $A$ is equivalent to a MP manifold $P\to A^*[1]$ together with a morphism $P\dasharrow\mathbf{G}_{small}$, i.e.\ to a morphism of Manin pairs
\begin{equation}\label{qpgb2}(\EE,A)\dasharrow(\dd,\on{diag}(\g)).\end{equation}

The groupoid multiplication defines a morphism of Manin pairs $$(\dd,\on{diag}(\g))\times (\dd,\on{diag}(\g))\dasharrow (\dd,\on{diag}(\g)),$$ and the fusion product of two quasi-Poisson $\g$-structures $(\TT M_i,TM_i)\dasharrow (\dd,\on{diag}(\g)$ (for $i=1,2$) corresponds to the composition $$(\TT M_1,TM_1)\times(\TT M_2, TM_2)\dasharrow (\dd,\on{diag}(\g))\times (\dd,\on{diag}(\g))\dasharrow (\dd,\on{diag}(\g)).$$ 

\begin{remark}
The morphism of Manin pairs \eqref{morphManinPairDesc1} is full if and only if the corresponding quasi-Poisson $\g$-structure on $M$ is non-degenerate \cite{Alekseev99,Alekseev00}, namely \eqref{nondegenqh} holds.
\end{remark}

\begin{remark}
The equivalence between morphisms of Manin pairs \eqref{morphManinPairDesc1} and quasi-Poisson $\g$ structures is just a restatement of the results in \cite{QuasiPoissonAs} (in particular, see the first paragraph of \cite[\S~5]{QuasiPoissonAs}).
\end{remark}

\begin{remark}[Special case of Theorem~\ref{mainthm1}] \label{establishedResult}

A referee explained to us that Theorem~\ref{mainthm1} was already established when there was a moment map $\Phi:M\to G$ for the quasi-Poisson $\g$-manifold $(M,\rho,\pi)$. 

In this case, \cite{Bursztyn03} describes an embedding of the Lie algebroid $T^*M$ as a Dirac structure $L\subset TM\oplus T^*M$  (see Remark~\ref{LieAlgL}). Let $s,t:\Gamma\to M$ be the source and target maps of the Lie groupoid integrating $T^*M\cong L$.
Then \cite{Bursztyn08}, there is a full morphism of Manin pairs $$(s\times t,K_0):(\mathbb{T}\Gamma,T\Gamma)\dasharrow (\mathbb{T}M\times \overline{\mathbb{T}M},L\times L).$$ This follows from \cite{Bursztyn03-1} or \cite{Ponte05} together with \cite{PureSpinorsOnL}.

 There are also full morphisms of Manin pairs $$(\Phi\times\Phi^{-1},K_1):(\mathbb{T}M\times \overline{\mathbb{T}M},L\times L)\dasharrow(\A_G\times \A_G, E_G\times E_G),$$ describing the Hamiltonian quasi-Poisson $\g$-structure on $M$ \cite{Bursztyn03,PureSpinorsOnL}, and $$(\on{mult},K_2):(\A_G\times \A_G, E_G\times E_G)\dasharrow (\A_G,E_G),$$ describing the fusion \cite{PureSpinorsOnL}.

The composition of these morphisms 
$$(\mathbb{T}\Gamma,T\Gamma)\dasharrow (\A_G,E_G)$$
describes a Hamiltonian quasi-symplectic $\g$-structure on $\Gamma$ \cite{Bursztyn03-1}.
\end{remark}

\subsection{Alternative proof of Theorem~\ref{mainthm1}}
As an application of Theorem~\ref{inffull}, we may sketch an alternative proof to Theorem~\ref{mainthm1}.

Let $(A,\rho,\DIF)$ be a quasi-Poisson $\g$-bialgebroid. Recall from Proposition~\ref{p-qpgb-htpa} that a quasi-Poisson $\g$-bialgebroid structure on $A$ is equivalent to a Poisson structure of degree $-1$ on $A[1]$ together with a Lie bialgebra action of $\hat\g$ on $A[1]$. However \cite{Xu95}, a Lie bialgebra action $\hat\g$ on $A[1]$ is equivalent to a morphism of Lie bialgebroids $T^*[1](A[1])\to \hat\g^*[1]$, Using the canonical symplectomorphism $T^*[1]A[1]\cong T^*[1]A^*$, we may rewrite this as \begin{equation}\label{noR2morph}T^*[1]A^*\to \hat\g^*[1].\end{equation} $\hat\g^*[1]$ is just the MP-algebra $\mathbb{R}[2]\oplus\g[1]\oplus\bar\g$ described in \S~\ref{qhqpgmr}, so \eqref{noR2morph} is canonically equivalent (as in \S~\ref{atqps}) to a morphism of MP-algebroids \begin{equation}\label{infmorph1}T^*[1]A^*\times\mathbb{R}[2]\dasharrow \hat\g^*[1].\end{equation}
Since $T^*[1]A^*\times\mathbb{R}[2]$ is just the MP-algebroid corresponding to the Manin pair $(\TT A^*,T A^*)$ and $\hat\g^*[1]$ corresponds to the Manin pair $(\TT \g^*,\on{Gr}_{\pi_\g})$, where $\pi_\g$ is the Kirillov bivector field and $\on{Gr}_{\pi_\g}$ is the graph of $\pi_\g^\sharp:T^*\g^*\to T\g^*$; \eqref{infmorph1} just corresponds to a morphism of Manin pairs \begin{equation}\label{infmorph2}(\TT A^*,T A^*)\dasharrow(\TT \g^*,\on{Gr}_{\pi_\g}).\end{equation} This is equivalent to a Poisson structure $\pi_{A^*}$ on $A^*$ such that $\mu_\rho:(A^*,\pi_{A^*})\to(\g^*,\pi_\g)$ is a Poisson morphism (i.e. a moment map). It is not difficult to check that $\pi_{A^*}$ is just the linear Poisson structure on $A^*$ corresponding to the Lie algebroid structure on $A$. \eqref{infmorph2} is full if and only if $(A^*,\pi_{A^*})$ is actually a symplectic manifold, which implies that $A\cong TM$ as Lie algebroids.

Suppose $\Gamma$ is a source 1-connected groupoid integrating the Lie algebroid $A^*$. Example~\ref{T*A} (with $B=A^*$) and the fact that $\hat\g^*[1]$ integrates to $\mathbf{G}_{big}$ show that \eqref{infmorph1} integrates to the morphism of MP-groupoids
\begin{equation}\label{globmorph1}T^*[1]\Gamma\times\mathbb{R}[2]\dasharrow\mathbf{G}_{big}\end{equation} describing the Hamiltonian quasi-Poisson $\g$-structure on $\Gamma$. In the language of Manin pairs, this is a morphism \begin{equation}\label{globmorph2}(\TT\Gamma,T\Gamma)\dasharrow(\A_G,E_G).\end{equation}

Theorem~\ref{inffull} states that \eqref{globmorph2} is full if and only if \eqref{infmorph2} is full. However \eqref{globmorph2} is full if and only if $\Gamma$ is a Hamiltonian quasi-symplectic $\g$-groupoid (see \S~\ref{atqps}), while \eqref{infmorph2} is full if and only if $A\cong TM$ as Lie algebroids. Consequently, in light of Remark~\ref{QPBAisQPM}, source 1-connected Hamiltonian quasi-symplectic $\g$-groupoids are in one-to-one correspondence with integrable quasi-Poisson manifolds. One can prove the rest of Theorem~\ref{mainthm1} by simply checking the details in the above argument.

\section{Hamiltonian quasi-Poisson $\g$-groups}\label{sec:fex}
Since quasi-Poisson $\g$-bialgebroids are equivalent to morphisms of Manin pairs \eqref{qpgb2}, it follows that 1-connected Hamiltonian quasi-Poisson $\g$-groups are classified by morphisms of Manin pairs
\begin{equation}(\id,K):\label{hqpg2}(\mf{f},\mf{h})\dasharrow(\dd,\on{diag}(\g)),\end{equation}
where $\mf{f}$ is a quadratic Lie algebra, $\h$ is a Lagrangian subalgebra, and $\id:\ast\to\ast$ is identity map for the point.

\begin{proposition}
Morphisms of Manin pairs of the form \eqref{hqpg2} are equivalent to quadruples $(\mf{f},\h\subset\mf{f},\h^*\subset\mf{f},\rho:\g\to \h)$, where $\mf{f}$ is a quadratic Lie algebra, and $\h,\h^*\subset \mf{f}$ are two subalgebras such that $\h$ is Lagrangian and $\mf{f}=\h\oplus \h^*$ as a vector space. Furthermore $\rho:\g\to \h$ is a Lie algebra morphism that satisfies
\begin{enumerate}
\item $\rho^*:\h^*\to \g$ is a Lie algebra morphism,
\item $\langle x,y\rangle_{\mf{f}}=\langle \rho^*(x),\rho^*(y)\rangle_\g$, and
\item $[\rho(\g),\h^*]\subset \h^*$ with $\rho^*[\rho(\xi),x]=[\xi,\rho^*(x)]$ for $\xi\in \g$ and $x\in \h^*$
\end{enumerate}
\end{proposition}

\begin{remark} It should be clear that these conditions define a Lie bialgebra morphism $\hat\g\to \hat\h$, where the Lie algebra bracket on $(\hat\h)^*[1]\cong \mathbb{R}[2]\oplus\h^*[1]\oplus\h^*$ is given on $\h^*[1]$ by the quadratic form on $\mf{f}$, and $\mathbb{R}[2]$ is central.

\end{remark}

\begin{proof}
Since $\dd=\g\oplus\bar\g$ is the vector space direct sum of the subalgebras $\on{diag}(\g)$ and $\g$, it follows that $\mf{f}$ is a vector space direct sum of the subalgebras $\h$ and 
$$\g\circ K=\{x\in \mf{f}\mid (y,x)\in K\subset \dd\oplus\bar{\mf{f}}, \text{ for some } y\in \g\}.$$
 We may identify $\g\circ K$ with $\h^*$ using the quadratic form on $\mf{f}$. Then $K\subset \dd\oplus\bar{\mf{f}}=\g+\bar\g+\overline{\h+\h^*}$ can be written as \begin{equation}\label{Kequat}K=\{(\xi,\xi,\rho(\xi),0)\in \g+\bar\g+\overline{\h+\h^*}\}+\{(\rho^*(x),0,0,x)\in\g+\bar\g+\overline{\h+\h^*}\},\end{equation} where $\rho:\g\to \h$ is a Lie algebra morphism.

On the other hand, suppose that $(\mf{f},\h\subset\mf{f},\h^*\subset\mf{f},\rho:\g\to \h)$ is a quadruple satisfying the assumptions. Then \eqref{Kequat} is a Lagrangian subalgebra of $\dd\oplus\bar{\mf{f}}$ defining a morphism of Manin pairs \eqref{hqpg2}.

\end{proof}

\appendix

\section{$\mathcal{VB}$-groupoids}\label{vbgroupoids}
\begin{definition}
A $\mathcal{VB}$-groupoid is a Lie groupoid in the category of smooth vector bundles (or a vector bundle in the category of Lie groupoids). In more detail, it is a diagram of the form
\begin{equation}\label{VBdiag}\xymatrix{
V\ar[r]^{\hat q}\ar@<-0.5ex>[d]_{\tilde s}\ar@<0.5ex>[d]^{\tilde t}& \Gamma\ar@<-0.5ex>[d]_s\ar@<0.5ex>[d]^t\\
E\ar[r]_q&M
}\end{equation}
where $\hat q:V\to\Gamma$ and $q:E\to M$ are vector bundles, and $V\rightrightarrows E$ and $\Gamma\rightrightarrows M$ are Lie groupoids whose source, target, multiplication and object inclusion maps ($\tilde s,\tilde t,\tilde m, \tilde i$ and $s,t,m,i$, respectively) are morphisms of vector bundles.

%

\end{definition}

Suppose $\Gamma\rightrightarrows M$ is a Lie groupoid. Then applying the tangent functor, we get a $\mathcal{VB}$-groupoid $T\Gamma\rightrightarrows TM$.  

\cite[Proposition 11.2.5]{Mackenzie05} states that if $V\to\Gamma$ is a $\mathcal{VB}$-groupoid, then $V^*\to\Gamma$ naturally inherits the structure of a $\mathcal{VB}$-groupoid. Briefly, if $(\on{Gr}_{\tilde m})_{(f,g,g)}$ is the fibre of the graph of the multiplication for $V$ at the point $(f,g,h)\in\on{Gr}_m$, then \begin{multline}(\on{Gr}_{\tilde m})^\perp_{(f,g,h)}=\{(\alpha,\beta,\gamma)\in V^*_f\times V^*_g\times V^*_h\mid 0=\alpha(u)+\beta(v)-\gamma(w),\\\forall (u,v,w)\in(\on{Gr}_{\tilde m})_{(f,g,h)}\}\end{multline}
is the fibre of the graph of the multiplication for $V^*$ at the point $(f,g,h)$. Meanwhile, if $g\in i(M)\subset \Gamma$ is an identity element, and $\tilde i(E_g)\subset V_g$ is the fibre of the identity elements of $V$ over $g$, then $$\tilde i(E_g)^\perp=\{\alpha\in V^*_g\mid \alpha(v)=0\;\;\forall v\in \tilde i(E_g)\}$$ is the fibre of the identity elements of $V^*$ over $g$. With this structure $V^*$ is called the dual $\mathcal{VB}$-groupoid.

Consequently $T^*\Gamma$ also has the structure of a $\mathcal{VB}$-groupoid.

\begin{remark}[Technical note]
The Theorems in \cite[\S~11.2]{Mackenzie05} assume that the ``double source condition'' is satisfied for the $\mathcal{VB}$-groupoids involved. That is to say if \eqref{VBdiag}
%
%
%
 is a $\mathcal{VB}$-groupoid, then the ``double source map'' \begin{equation}\label{dbsource}(\hat q,\tilde s):V\to \Gamma\times_{s,q} E\end{equation} is a surjective submersion. In order to apply these theorems to the $\mathcal{VB}$-groupoids used in our paper, we need the following lemma.

\begin{lemma}
If \eqref{VBdiag} is a Lie groupoid in the category of smooth vector bundles, then \eqref{dbsource} is a surjective submersion.
\end{lemma}

\begin{proof}
We begin by showing that \eqref{dbsource} is surjective. View $\Gamma$ as the zero section of the vector bundle $V\to \Gamma$, and let $g\in \Gamma$. The vector space $T_gV$ decomposes into directions tangent to the fibres and directions tangent to the zero section, namely
$$T_gV=V_g\oplus T_g\Gamma.$$ Similarly $T_{s(g)}E$ has a natural decomposition $$T_{s(g)}E=E_{s(g)}\oplus T_{s(g)}M.$$ Since $\tilde s$ is a morphism of vector bundles, $T_g\tilde s$ decomposes as the direct sum $$T_g\tilde s=\tilde s\rvert_g\oplus T_g s:V_g\oplus T_g\Gamma\to E_{s(g)}\oplus T_{s(g)}M.$$ However $V$ was assumed to be a Lie groupoid, hence $\tilde s$ is a surjective submersion, and consequently $\tilde s\rvert_g:V_g\to E_{s(g)}$ is surjective. It follows that \eqref{dbsource} is surjective.

To show that \eqref{dbsource} is also a submersion, apply the previous argument to
$$\xymatrix{
TV\ar[r]^{T\hat q}\ar@<0.5ex>[d]^{T\tilde s}\ar@<-0.5ex>[d]_{T\tilde t}& T\Gamma\ar@<0.5ex>[d]^{Ts}\ar@<-0.5ex>[d]_{Tt}\\
TE\ar[r]^{Tq}&TM
}$$
\end{proof}
\end{remark}


\bibliography{basicbib}{}
\bibliographystyle{hplain2}

\end{document}